\documentclass[MSNbibl,number,citesort,dvips]{arxbj}
\usepackage{upgreek}
\usepackage{graphicx}

\aid{0}
\volume{19}
\issue{5B}
\pubyear{2013}
\firstpage{2359}
\lastpage{2388}
\doi{10.3150/12-BEJ455} %kopijuoti is PTS

\makeatletter
\def\eqref#1{(\ref{#1})}
\newcommand{\rrVert}{\Vert}
\newcommand{\rrvert}{\vert}
\newcommand{\llVert}{\Vert}
\newcommand{\llvert}{\vert}
\newtheorem{theorem}{Theorem}%[section]
\newtheorem{proposition}[theorem]{Proposition}
\newtheorem{lemma}{Lemma}
\newremark{remark}{Remark}
\makeatother

\begin{document}
\begin{frontmatter}

\title{Optimal hypothesis testing for high dimensional covariance matrices}
\runtitle{Testing high dimensional covariance matrices}

\begin{aug}
%%%% inicialai - be tarpu
\author{\fnms{T.~Tony} \snm{Cai}\thanksref{e1}\ead[label=e1,mark]{tcai@wharton.upenn.edu}} \and
\author{\fnms{Zongming} \snm{Ma}\corref{}\thanksref{e2}\ead[label=e2,mark]{zongming@wharton.upenn.edu}}
\runauthor{T.T. Cai and Z. Ma} %% auto
\address{Department of Statistics, The Wharton School, University of
Pennsylvania, Philadelphia, PA 19104, USA. \printead{e1,e2}}
\end{aug}

% HISTORY:
\received{\smonth{10} \syear{2011}}
\revised{\smonth{5} \syear{2012}}

% ABSTRACT
%
\begin{abstract}
This paper considers testing a covariance matrix $\Sigma$ in the high
dimensional setting where the dimension $p$
can be comparable or much larger than the sample size $n$. The problem
of testing the hypothesis $H_0\dvtx\Sigma= \Sigma_0$
for a given covariance matrix $\Sigma_0$ is studied from a minimax
point of view. We first characterize the boundary that separates the
testable region from the non-testable region by the Frobenius norm when
the ratio between the dimension $p$ over the sample size $n$ is
bounded. A test based on a $U$-statistic is introduced and is shown to
be rate optimal over this asymptotic regime. Furthermore, it is shown
that the power of this test uniformly dominates that of the corrected
likelihood ratio test (CLRT) over the entire asymptotic regime under
which the CLRT is applicable. The power of the $U$-statistic based test
is also analyzed when $p/n$ is unbounded.
\end{abstract}

% KEYWORDS
%

\end{frontmatter}
%

%s1 ###
\section{Introduction}
\label{sec:intro}

%Many statistical procedures in multivariate analysis require
%structural assumptions on the covariance matrix.

Covariance structure plays a fundamental role in multivariate analysis
and testing the covariance matrix is an important problem. Let
$X_1,\dots, X_n$ be $n$ independent and identically distributed
$p$-vectors following a multivariate normal distribution $N_p(0,\Sigma
)$. A~hypothesis testing problem of significant interest is testing
\begin{equation}
\label{eq:H0-cov} H_0\dvt \Sigma= I. % \mbox{vs.}  H_1: \covmat\in\Theta,
\end{equation}
Note that any null hypothesis $H_0\dvt \Sigma= \Sigma_0$ with a given
positive definite covariance matrix $\Sigma_0$ is equivalent to
\eqref{eq:H0-cov}, since one can always transform $X_i$ to $\tilde
{X}_i = \Sigma_0^{-1/2}X_i$ and then test \eqref{eq:H0-cov} based on
the transformed data.
% \nb{remark on $H_0: \covmat=\covmat_0$. }

This testing problem has been well studied in the classical setting of
small $p$ and large $n$. See, for example, Anderson \cite{ande03} and
Muirhead \cite{muir82}. In particular, the likelihood ratio test (LRT)
is commonly used. Driven by a wide range of contemporary scientific
applications, analysis of high dimensional data is of significant
current interest. In the high dimensional setting, where the dimension
can be comparable to or even much larger than the sample size, the
conventional testing procedures such as the LRT perform poorly or are
not even well defined. Several testing procedures designed for the
high-dimensional setting have been proposed. Let $S = \frac{1}{n}\sum_{i=1}^n X_i X_i'$ be the sample covariance matrix.
The existing tests for \eqref{eq:H0-cov} in the literature can be
categorized as the following according to the asymptotic regime under
which they are suitable:
\begin{itemize}
\item$p$ fixed and $n\to\infty$. In this classical asymptotic
regime, conventional tests for \eqref{eq:H0-cov} include the
likelihood ratio test (LRT) \cite{ande03}, Roy's largest root test
\cite{roy57}, and Nagao's test \cite{naga73}. In particular, the LRT
statistic is $\mathit{LR}_n = n L_n$, where
\[
L_n = \operatorname{tr} {S}-\log\operatorname{det}(S)-p.
\]
The asymptotic distribution of $\mathit{LR}_n$ under\vspace*{1pt} $H_0$ is $\chi^2_{p(p+1)/2}$.

% \item Both $n,p\to\infty$ and $p/n\to0$.
% The case is not as well studied as the classical setting of fixed $p$
%and large $n$. Srivastava \cite{sriv05} proposed a test that is valid
%for the case when $n = O(p^\delta)$, $0<\delta\leq1$, which includes
%the case $p/n\to0$. {\bf\red$\leftarrow$ Why this is true? Any other
%paper on this case?} {\bf\blue Alternatively, we can say the
%following: This regime is not well studied specifically. However,
%several tests proposed for the high dimensional setting are also
%suitable for this case. The tests include \ldots } {\bf\red I believe no
%one has specifically studied this case. So another option is simply to
%delete the case from this list.}

\item Both $n,p\to\infty$ and $p/n\to c\in(0,\infty)$.
Investigation in this asymptotic regime has been very active in the
past decade. For example, Johnstone \cite{john01} revisited Roy's
largest root test and derived the Tracy--Widom limit of its null
distribution. Ledoit and Wolf \cite{lewo02} proposed a new test based
on Nagao's proposal. See also Srivastava \cite{sriv05}. When $p$
grows, the chi-squared limiting null distribution of the LRT statistic
$\mathit{LR}_n$ is no longer valid. Recently, Bai \textit{et al.} \cite{baji09} proposed
a corrected LRT when $c < 1$, and Jiang \textit{et al.} \cite{jiji10} extended
it to the case when $p<n$ and $c = 1$. Here, for $c_n = p/n$, the test
statistic of the corrected LRT is
\begin{equation}
\label{eq:clrt} \mathit{CLR}_n = \frac{L_n - p[1-(1-c_n^{-1}) \log(1-c_n)] - {1}/{2}\log
(1-c_n)}{\sqrt{-2\log(1-c_n)-2c_n}},
\end{equation}
whose asymptotic null distribution is $N(0,1)$. Note that no test based
on the likelihood ratio can be defined when $p>n$ or $c > 1$.

\item Both $n,p\to\infty$ and $p/n\to\infty$.
% \nb{Need to be revised.}
This is the ultra high-dimensional setting and both the LRT and
corrected LRT are not well defined in this case. The testing problem in
this asymptotic regime is not as well studied as in the previous
categories. % Srivastava \cite{sriv05} proposed a test that is valid
%for the case when $n = O(p^\delta)$, $0<\delta\leq1$, which includes
%the case $p/n\to0$. {\bf\red$\leftarrow$ Double check!}
Birke and Dette \cite{bide05} derived the asymptotic null distribution
of the Ledoit--Wolf test under the current asymptotic regime. More
recently, Chen \textit{et al.} \cite{chzh10} proposed a new test statistic and
derived its asymptotic null distribution when both $n,p\to\infty$,
regardless of the limiting behavior of $p/n$.
\end{itemize}

When the dimension $p$ grows together with the sample size $n$, the
focus of most of the aforementioned papers is mainly on finding the
asymptotic null distribution of the proposed test statistic, so the
significance level of the test can be controlled. The few exceptions
include Srivastava \cite{sriv05} and Chen \textit{et al.} \cite{chzh10}, where
the asymptotic pointwise power of the proposed tests is also studied.
Recently, Onatski \textit{et al.} \cite{onmo11} established the regime of
mutual contiguity of the joint distributions of the sample eigenvalues
under the null and under the special alternative of rank one
perturbation to the identity matrix, and then applied Le Cam's third
lemma to study the pointwise power of a collection of eigenvalue based
tests for \eqref{eq:H0-cov} against this special class of alternative.

In the present paper, we investigate this testing problem in the
high-dimensional settings from a minimax point of view. Consider
testing \eqref{eq:H0-cov} against a composite alternative hypothesis
\begin{equation}
\label{eq:H1-cov} H_1\dvt \Sigma\in\Theta,\qquad  \mbox{where }
\Theta= \Theta_n = \bigl\{\Sigma\dvt \llVert \Sigma-I
\rrVert_F\geq\varepsilon_n \bigr\}.
\end{equation}
Here, $\llVert  A \rrVert_F=(\sum_{ij} a_{ij}^2)^{1/2}$ denotes the Frobenius
norm of a matrix $A=(a_{ij})$. It is clear that the difficulty of
testing between $H_0$ and $H_1$ depends on the value of $\varepsilon_n$: the
smaller $\varepsilon_n$ is, the harder it is to distinguish between
the two
hypotheses. An interesting question is: What is the boundary that
separates the testable region, where it is possible to reliably detect
the alternative based on the observations, from the untestable region,
where it is impossible to do so? This problem is connected to the
classical contiguity theory. It is also important to construct a test
that can optimally distinguish between the two hypotheses in the
testable region. The high-dimensional settings here include all the
cases where the dimension $p=p_n\to\infty$ as the sample size $n \to
\infty$, and there is no restriction on the limit of $p/n$ unless
otherwise stated.

For a given the significance level $0<\alpha<1$, our first goal is to
identify the separation rate $\varepsilon_n$ at which there exists a test
$\phi$ based on the random sample $\{X_1, \ldots , X_n\}$ such that
\[
\inf_{\Sigma\in\Theta}\mathsf{P}_{\Sigma}(\phi\mbox{ rejects }
H_0) \geq\beta> \alpha.
\]
Hence, the test is able to detect any alternative that is separated
away from the null by a certain distance $\varepsilon_n$ with a guaranteed
power $\beta> \alpha$. Our second goal is to construct such a testing
procedure $\phi$.

The major contribution of the current paper is threefold. First, we
show that if $p/n$ is bounded, then the rate $\varepsilon_n$ needs to
be no
less than $b\sqrt{p/n}$ for some constant $b$. In addition, it is
shown that if $\varepsilon_n = b\sqrt{p/n}$, there exists a test
$\psi$ of
significance level $\alpha$, such that $\lim_{n\to\infty}\inf_{\Theta}\mathsf{P}_{\Sigma}(\mbox{$\psi$ rejects $H_0$}) >\alpha
$, and
the power tends to $1$ if $b = b_n\to\infty$. The test is motivated
by the proposal in Chen \textit{et al.} \cite{chzh10}. (We use $\psi$ to
denote the specific test that we construct, while $\phi$ is used to
denote a generic test.)
Here, we no longer require $p/n$ to be bounded, and the explicit
expression for the asymptotic power of $\psi$ is also given. Moreover,
we show that the asymptotic power of $\psi$ on $\Theta_n$ uniformly
dominates that of the corrected LRT by Bai \textit{et al.} \cite{baji09} and
Jiang \textit{et al.} \cite{jiji10} over the entire asymptotic regime under
which the corrected LRT is defined, that is, $p<n$ and $p/n\to c\in(0,1]$.

The rest of the paper is organized as the following. In Section \ref
{sec:lowbd}, after introducing basic notation and definitions, we
establish a lower bound of the separation rate $\varepsilon_n$.
Section \ref
{sec:uppbd} introduces the test based on a $U$-statistic and provides a
Berry--Essen bound for its weak convergence to the normal limit under
both the null and the alternative hypotheses, which leads to the
establishment of its guaranteed power over $\Theta$ when $\varepsilon_n =
b\sqrt{p/n}$. Furthermore, we also show that the power of this test
uniformly dominates that of the corrected LRT. The theoretical results
are supported by the numerical experiments in Section \ref{sec:numeric}.
Further discussions on the connections of our results and those of
related testing problems are given in Section \ref{sec:disc}. The main
results are proved in Section~\ref{sec:proof}.

%s2 ###
\section{Lower bound}
\label{sec:lowbd}

In this section, we establish a lower bound for the separation rate
$\varepsilon_n$ in \eqref{eq:H1-cov}. The result in Section~\ref{sec:uppbd}
will show that this lower bound is rate-optimal. The lower and upper
bounds together characterize the separation boundary between the
testable and non-testable regions when the ratio of the dimension $p$
over the sample size $n$ is bounded. This separation boundary can then
be used as a minimax benchmark for the evaluation of the performance of
a test in this asymptotic regime.

% Recall that we are interested in the high-dimensional setting where
%the dimension $p=p_n\to\infty$ as the sample size $n \to\infty$. So,
%the dimension $p_n$, the covariance matrix $\covmat= \covmat_n$, the
%alternative class $\Theta= \Theta_n$ and many other quantities are
%all indexed by $n$. The subscript will be shown only when clarity
%dictates.
We begin with basic notation and definitions. Throughout the paper,
a test $\phi= \phi_n(X_1,\dots, X_n)$ refers to a measurable
function which maps $X_1,\dots, X_n$ to the closed interval $[0,1]$,
where the value stands for the probability of rejecting $H_0$. So, the
significance level of $\phi$ is $\mathsf{P}_{I}(\mbox{$\phi$ rejects
$H_0$}) = \mathsf{E}_{I}\phi$, and its power at a certain alternative
$\Sigma$ is $\mathsf{P}_{\Sigma}(\mbox{$\phi$ rejects $H_0$}) =
\mathsf{E}_{\Sigma}\phi$. Here and after, $\mathsf{P}_{\Sigma}, \mathsf
{E}_{\Sigma},
\operatorname{\mathsf{Var}}_{\Sigma}$ and $\operatorname{\mathsf
{Cov}}_{\Sigma}$ denote the induced probability
measure, expectation, variance and covariance when $X_1,\dots, X_n
\stackrel{\mathrm{i.i.d.}}{\sim}N_p(0,\Sigma)$.
The subscript is shown only when clarity dictates.
% When there is no ambiguity, we sometimes drop the subscript $\covmat$.

To state the lower bound result, let $\varepsilon_n = b\sqrt{p/n}$
for some
constant $b$, and define
\begin{equation}
\label{eq:class-lowbd} % \Theta=
\Theta(b) = \bigl\{\Sigma\dvt \llVert \Sigma-I
\rrVert_F \geq b\sqrt {p/n} \bigr\}.
\end{equation}

\begin{theorem}[(Lower bound)]
\label{thm:lowbd}
Let $0<\alpha< \beta< 1$. Suppose that as $n\to\infty$, $p\to
\infty$ and that $p/n\leq\kappa$ for some constant $\kappa<\infty$
and all $n$. Then there exists a constant $b = b(\kappa,\beta-\alpha
) < 1$, such that for any test $\phi$ with significance level $\alpha
$ for testing $H_0\dvt \Sigma= I$,
\[
\limsup_{n\to\infty}\inf_{\Sigma\in\Theta(b)} \mathsf {E}_{\Sigma
}\phi< \beta.
\]
\end{theorem}

Theorem \ref{thm:lowbd} shows that no level $\alpha$ test for \eqref
{eq:H0-cov} can distinguish between the two hypotheses with power
tending to $1$ as $n$ and $p$ grow, when the separation rate
$\varepsilon_n$
is of order $\sqrt{p/n}$. Hence, it provides a lower bound for the
separation rate.

We now give an outline of the proof for Theorem \ref{thm:lowbd}, while
the complete proof is provided in Section \ref{sec:pf-lowbd}. Consider
the following ``least favorable'' subset of $\Theta(b)$:
\begin{equation}
\label{eq:least-favor-class} \Theta^*(b) = \biggl\{\Sigma_{v} =
\biggl[1-\frac{b}{\sqrt {n(p-1)}} \biggr] I_{p\times p} + \frac{b}{\sqrt{n(p-1)}}vv'\dvt
v\in\{ \pm1\}^p \biggr\}.
\end{equation}
With slight abuse of notation,
let $P_0$ be the probability measure when $X_1,\dots X_n\stackrel
{\mathrm{i.i.d.}}{\sim}N_p(0,I)$
and $P_v$ the probability measure when $X_1,\dots, X_n\stackrel
{\mathrm{i.i.d.}}{\sim}
N_p(0,\Sigma_v)$. In addition, let $P_1 = \frac{1}{2^p}\sum_{v\in\{
\pm1\}^p} P_v$ be the average measure of the $P_v$'s. Then for any
test $\phi$, the sum of probabilities of its two types of errors satisfies
\begin{eqnarray*}
\sup_{v}\mathsf{E}_0\phi+ \mathsf{E}_v(1-
\phi) & \geq&\inf_\psi\sup_{v} \mathsf{E}_0\psi+
\mathsf{E}_v(1-\psi)
\\
& \geq&\inf_\psi\frac{1}{2^p}\sum_{v}
\mathsf{E}_0\psi+ \mathsf {E}_v(1-\psi)
\\
& =& \inf_\psi\mathsf{E}_0\psi+ \mathsf{E}_1(1-
\psi)
\\
& =& 1 - \frac{1}{2}\|P_1 - P_0\|_1.
\end{eqnarray*}
Here, $\mathsf{E}_0,\mathsf{E}_v$ and $\mathsf{E}_1$ denote the
expectation under $P_0$,
$P_v$ and $P_1$ respectively, and $\|P_1-P_0\|_1$ is the $L_1$ distance
between $P_0$ and $P_1$.
Thus, we obtain
\[
\inf_{\Sigma\in\Theta(b)} \mathsf{E}_{\Sigma}\phi\leq\inf_v
\mathsf{E}_v\phi \leq\mathsf{E}_0\phi+ \frac{1}{2}
\|P_1-P_0\|_1 = \alpha+ \frac
{1}{2}\|
P_1-P_0\|_1.
\]
To control the rightmost side, we bound the $L_1$ distance by the
chi-square divergence as
\[
\|P_1 - P_0\|_1^2 \leq
\mathsf{E}_0\biggl\llvert \frac{\mathrm{d} P_1}{\mathrm{d} P_0}-1\biggr
\rrvert^2 = \mathsf{E}_0\biggl\llvert \frac{\mathrm{d} P_1}{\mathrm{d} P_0}
\biggr\rrvert^2 - 1 = \int\frac
{f^2_1}{f_0} - 1,
\]
where $f_i$ is the density function of $P_i$ for $i=0,1$.
So, the proof can be completed by showing that for an appropriate
choice of the constant $b$, one obtains $\int\frac{f_1^2}{f_0} -
1\leq4(\beta-\alpha)^2$.

%We conclude this section with two remarks:

\begin{remark}
\label{correlation.remark}
(a). Note that all the covariance matrices in the least favorable
configuration $\Theta^*(b)$ defined in \eqref{eq:least-favor-class}
have diagonal elements all equal to $1$. Thus, they are also
correlation matrices. So the proof of Theorem \ref{thm:lowbd} readily
establishes an analogous lower bound result on testing $H_0\dvt R = I$
with $R$ the population correlation matrix.
% For details, see Theorem \ref{thm:lowbd-corr} in Section

(b). The lower bound argument here does not extend to the case when
$p/n$ is unbounded, because the chi-square divergence becomes unbounded.
\end{remark}

%s3 ###
\section{Upper bound}
\label{sec:uppbd}

% Theorem \ref{thm:lowbd} showed that $H_0$ cannot be tested perfectly
%against the composite alternative in \eqref{eq:class-lowbd}, and so $
%in terms of the Frobenius distance.
In this section, we show that there exists a level $\alpha$ test whose power
over $\Theta_n$ is uniformly larger than a prescribed value $\beta>
\alpha$, if $\varepsilon_n = b\sqrt{p/n}$ for a large enough
constant $b$.
This matches the lower bound result in Theorem \ref{thm:lowbd} when
$p/n$ is bounded. In addition, the results in the current section
remain valid even when $p/n$ is unbounded.
% satisfying $\fnorm{\covmat- I} \geq b\sqrt{p/n}$ for some $b>0$.

We first introduce the test statistic in Section \ref{sec:chen-test},
followed by a study on the rate of convergence of its distribution to
the normal limit under both the null and the alternative hypotheses in
Section~\ref{sec:rate}. Section \ref{sec:power} then uses the rate of
convergence result to study the asymptotic power of the proposed test.
Finally, Section \ref{sec:CLRT} shows that the test dominates the
corrected LRT in \eqref{eq:clrt} when $p/n\to c\in(0,1]$.

%s3.1 ###
\subsection{Test statistic}
\label{sec:chen-test}

Given a random sample $X_1,\dots, X_n\stackrel{\mathrm{i.i.d.}}{\sim
}N_p(0,\Sigma)$, a natural
approach to test between \eqref{eq:H0-cov} and \eqref{eq:H1-cov} is
to first estimate the squared Frobenius norm $\llVert  \Sigma-I \rrVert_F^2 =
\operatorname{tr}(\Sigma-I)^2$ by some statistic $T_n = T_n(X_1,\dots
, X_n)$, and
then reject the null hypothesis if $T_n$ is too large. To estimate
$\llVert  \Sigma-I \rrVert_F^2 = \operatorname{tr}(\Sigma-I)^2$,
note that $\mathsf{E}_{\Sigma
} h(X_1,X_2) = \operatorname{tr}(\Sigma-I)^2$ where
\begin{equation}
\label{eq:h} h(X_1,X_2) = \bigl(X_1'X_2
\bigr)^2 - \bigl(X_1'X_1 +
X_2'X_2\bigr) + p.
\end{equation}
Therefore, $\operatorname{tr}(\Sigma-I)^2$ can be estimated by the
following $U$-statistic
\begin{equation}
\label{eq:T_n} T_n = \frac{2}{n(n-1)} \sum
_{1\leq i < j \leq n} h(X_i, X_j),
\end{equation}
for which we have
\begin{eqnarray}\label {eq:E-T_n}
\mu_n(\Sigma) & = &\mathsf{E}_{\Sigma}(T_n) =
\operatorname {tr}(\Sigma-I)^2,
\\\label{eq:var-T_n}
\sigma_n^2(\Sigma) & =& \operatorname{
\mathsf{Var}}_{\Sigma}(T_n) = \frac{4}{n(n-1)} \bigl[
\operatorname{tr}^2\bigl(\Sigma^2\bigr) +
\operatorname{tr}\bigl(\Sigma^4\bigr) \bigr] + \frac{8}{n}
\operatorname{tr} \bigl(\Sigma^2(\Sigma- I)^2 \bigr).
\end{eqnarray}
Here, verifying \eqref{eq:E-T_n} is straightforward, and \eqref
{eq:var-T_n} is proved in Appendix \ref{app:var-T-n}.
For the $U$-statistic $T_n$, the proof for Theorem 2 of Chen \textit{et al.}
\cite{chzh10} essentially established the following.

\begin{proposition}[(Theorem 2 of \cite{chzh10})]
\label{prop:T-n-CLT}
Suppose that $p\to\infty$ as $n\to\infty$. If a sequence of
covariance matrices satisfy $\operatorname{tr}(\Sigma^2)\to\infty$
and $\operatorname{tr}
(\Sigma^4)/\operatorname{tr}^2(\Sigma^2)\to0$ as $n\to\infty$,
then under
$\mathsf{P}_{\Sigma}$, we have
\[
\frac{T_n - \mu_n(\Sigma)}{\sigma_n(\Sigma)} \Rightarrow N(0,1).
\]
\end{proposition}

Note that as $p\to\infty$, the identity matrix $I_{p\times p}$
satisfies the condition of the above proposition. Also note that $\mu_n(I) = 0$ and $\sigma_n^2(I) = {4p(p+1)\over n(n-1)}$. Thus,
Proposition \ref{prop:T-n-CLT} quantifies the behavior of $T_n$ under
$H_0$, and we could define the test as the following: For any $\alpha
\in(0,1)$, an asymptotic level $\alpha$ test based on $T_n$ is given by
\begin{equation}
\label{eq:chen-test} \psi= I \biggl(T_n > z_{1-\alpha}\cdot2
\sqrt{{p(p+1)\over n(n-1)}} \biggr).
\end{equation}
%
% {\bf\red Is $\sigma_n(I)$ exactly equal to $\sigma_n^2(\covmat)$ when
%$\Sigma= I$? If so, it might be helpful to write it in a simpler
%expression, $\sigma^2(I) = {4p(p+1)\over n(n-1)}$.}
Here, $I(\cdot)$ is the indicator function, and $z_{1-\alpha}$
denotes the $100\times(1-\alpha)$th percentile of the standard normal
distribution. This test is motivated by the test introduced in Chen \textit{et al.} \cite{chzh10}, while the original proposal in \cite{chzh10}
involves higher order symmetric functions of the $X_i$'s.

In addition to specifying the rejection region in \eqref
{eq:chen-test}, Proposition \ref{prop:T-n-CLT} can also be used to
study the asymptotic power of $\psi$ over a sequence of simple
alternatives. However, to understand the power of $\psi$ over the
composite alternative $\Theta$ in \eqref{eq:H1-cov}, it is necessary
to understand the rate of convergence of $[T_n- \mu_n(\Sigma
)]/\sigma_n(\Sigma)$ to the normal limit, which is the central topic
of the next subsection.

% The test statistic we consider is adapted from Chen, et al
% \begin{align}
% \label{eq:V-n}
% V_n = \frac{1}{2}\sqrt{\frac{n(n-1)}{p(p+1)}} \bigg[ \frac{2}{n(n-1)}
%p \bigg].
% \end{align}
% For any $\alpha\in(0,1)$, a test of significance level $\alpha$
%based on $V_n$ is given by
% \begin{align}
% \label{eq:chen-test}
% \psi= I(V_n > z_{1-\alpha}).
% \end{align}
% Here, $I(\cdot)$ is the indicator function, and $z_{1-\alpha}$
%denotes the $100\times(1-\alpha)$th percentile of the standard normal
%distribution.
%
% To understand the behavior of $V_n$ and hence of $\psi$, define
% \begin{equation}
% \label{eq:test-cov}
% \begin{aligned}
% % T_n =
% T_n(\covmat) & = \frac{2}{n(n-1)}\sum_{1\leq i<j\leq n}(X_i'X_j)^2 -
% & = \frac{2}{n(n-1)}\sum_{1\leq i<j\leq n}h(X_i,X_j),
% \end{aligned}
% \end{equation}
% where $h(X_1, X_2) = (X_1'X_2)^2 - X_1'X_1 - X_2'X_2 - \tr{(
% % \begin{align}
% % \label{eq:h}
% % h(X_1, X_2) = (X_1'X_2)^2 - X_1'X_1 - X_2'X_2 - \tr{(\covmat^2)} + 2
% % \end{align}
% % So, $T_n$ is a $U$-statistic.
% It is straightforward to verify that $\Ex_{\covmat}[T_n(\covmat) ] =
%0$, because each individual $h(X_i,X_j)$ has mean $0$. We can also
%show that
% % \nb{seems not to agree with \citep{chzh10} on the second term!}
% \begin{align}
% \label{eq:var-T_n}
% % \sigma_n^2 =
% \sigma_n^2(\covmat) = \var_{\covmat}[T_n(\covmat)] = \frac{4}{n(n-1)}
% \end{align}
% With $T_n(\covmat)$ and $\sigma_n^2(\covmat)$, $V_n$ can be expressed
%as
% \begin{align}
% \label{eq:V-n-T-n}
% V_n = T_n(I)/\sigma_n(I).
% \end{align}
%
% For $T_n(\covmat)$ and hence $V_n$,

%s3.2 ###
\subsection{Rate of convergence}
\label{sec:rate}

We now study the rate of convergence for the distribution of $[T_n-\mu_n(\Sigma)]/\sigma_n(\Sigma)$ to its normal limit in Kolmogorov distance.
% , i.e.,
% $\sup_{x\in\RR}\big| \Pr_{\covmat}(T_n(\covmat)/\sigma_n(\covmat)
Let $\Phi(\cdot)$ be the cumulative distribution function of the
standard normal distribution. We have the following Berry--Essen type bound.

\begin{proposition}
\label{prop:T-n-MGCLT}
Under the condition of Proposition \ref{prop:T-n-CLT}, there exists a
numeric constant $C$ such that
\[
\sup_{x\in\mathbb{R}}\biggl\llvert \mathsf{P}_{\Sigma} \biggl(
\frac{T_n
-\mu_n(\Sigma)}{\sigma_n(\Sigma)} \leq x \biggr) - \Phi(x) \biggr\rrvert \leq C \biggl[
\frac{1}{n} + \frac{\operatorname{tr}(\Sigma^4)}{\operatorname{tr}^2(\Sigma^2)} \biggr]^{1/5}.
\]
\end{proposition}

We outline the proof of Proposition \ref{prop:T-n-MGCLT} below, while
the complete proof is deferred to Section~\ref{sec:pf-MGCLT}.
The primary tool used in the proof is a Berry--Esseen type bound for
martingale central limit theorem by Heyde and Brown \cite{hebr70}.

We begin by giving a martingale representation of $T_n-\mu_n(\Sigma)$.
% For notational convenience, in what follows, $T_n$ and $\sigma_n^2$
%are used to abbreviate $T_n(\covmat)$ and $\sigma_n^2(\covmat)$.
Let $X_i\stackrel{\mathrm{i.i.d.}}{\sim}N_p(0,\Sigma)$.
Define filtration
\[
\mathcal{F}_0 = \sigma(\varnothing),\qquad   \mathcal{F}_k
= \sigma (X_1,\dots, X_k), \qquad  k=1,\dots, n.
\]
Also introducing the notation $\mathsf{E}_k[ \cdot ] = \mathsf
{E}_{\Sigma}[
\cdot  | \mathcal{F}_k]$. Then %$T_n$ can be represented as
\begin{equation}
\label{eq:T_n-MG} T_n - \mu_n(\Sigma) = \sum
_{k=1}^n \mathsf{E}_k
[T_n] - \mathsf {E}_{k-1}[T_n] = \sum
_{k=1}^n D_{nk}.
\end{equation}
Here, $\{D_{nk}\dvt k=1,\dots, n\}$ is a martingale difference sequence.
The explicit expression for $D_{nk}$ is
\begin{eqnarray}
\label{eq:T_n-MG-diff} D_{nk} &=& \frac{2}{n(n-1)} \bigl[
X_k'Q_{k-1}X_k - \operatorname
{tr}(Q_{k-1}\Sigma) \bigr]
\nonumber\\[-8pt]\\[-8pt]
&&{}  + \frac{2}{n} \bigl[
X_k'\Sigma X_k - \operatorname{tr}\bigl(
\Sigma^2\bigr) \bigr] - \frac{2}{n} \bigl[X_k'X_k
- \operatorname{tr}(\Sigma) \bigr],\nonumber  % v_{nk} - 2 u_{nk}.
% \frac{2}{n(n-1)}\Big[ X_k'Q_{k-1}X_k - \tr(Q_{k-1}\covmat) \Big] +
\end{eqnarray}
with $Q_{k-1} = \sum_{i=1}^{k-1}(X_iX_i' - \Sigma)$.
Let $\sigma_{nk}^2 = \mathsf{E}_{k-1}[D_{nk}^2]$, and we have $\sigma_n^2(\Sigma) = \sum_{k=1}^n \mathsf{E}_\Sigma[\sigma_{nk}^2]$.

Under the current setup, the main theorem in \cite{hebr70} specializes
to the following lemma.

\begin{lemma}
\label{lemma:MGCLT}
% Let $\sigma_{nk}^2 = \Ex_{k-1}[D_{nk}^2]$ be the associated
%increasing sequence. Recall that $\sigma_n^2(\covmat) = \var_
There exist a numeric constant $C$, such that
\begin{eqnarray}
\label{eq:mg-clt-bd} %
& & \sup_{x\in\mathbb{R}}\biggl
\llvert \mathsf{P}_{\Sigma} \biggl(\frac
{T_n-\mu_n(\Sigma)}{\sigma_n(\Sigma)} \leq x \biggr) -
\Phi(x) \biggr\rrvert
\nonumber\\[-8pt]\\[-8pt]
&&\quad \leq  C  \Biggl[ \frac{1}{\sigma_{n}^4(\Sigma)} \Biggl( \sum
_{k=1}^n\mathsf{E}_{\Sigma}
\bigl[D_{nk}^{4}\bigr] + \mathsf{E}_{\Sigma} \Biggl[
\sum_{k=1}^n\sigma_{nk}^2
- \sigma_n^2(\Sigma) \Biggr]^{2} \Biggr)
\Biggr]^{1/5}.\nonumber
\end{eqnarray}
\end{lemma}

Define
\begin{equation}
\label{eq:E1-E2} E_1 = \sum_{k=1}^n
\mathsf{E}_\Sigma\bigl[D_{nk}^4\bigr]
\quad \mbox{and}\quad   E_2 = \mathsf{E}_\Sigma \Biggl[\sum
_{k=1}^n \sigma_{nk}^2 -
\sigma_n^2(\Sigma) \Biggr]^{2}.
\end{equation}
The proof of Proposition \ref{prop:T-n-MGCLT} could then be completed
by showing that $E_1/\sigma_n^4(\Sigma) = \mathrm{O}(1/n)$ and $E_2/\sigma_n^4(\Sigma) = \mathrm{O}(\operatorname{tr}(\Sigma^4)/\operatorname
{tr}^2(\Sigma^2))$. See Section
\ref{sec:pf-MGCLT} for details.
% \nb{double check all the notation here!}

%s3.3 ###
\subsection{Power of the test}
\label{sec:power}
% \nb{consider changing section title!}
Equipped with Proposition \ref{prop:T-n-MGCLT}, we now investigate the
power of the test $\psi$ in \eqref{eq:chen-test} over the composite
alternative $H_1\dvt\Sigma\in\Theta(b)$, with $b<1$, where $\Theta
(b)$ is defined in \eqref{eq:class-lowbd}. In particular, we have the
following result.

\begin{theorem}[(Upper bound)]
\label{thm:uppbd}
Suppose that $p\to\infty$ as $n\to\infty$. For any significance
level \mbox{$\alpha\in(0,1)$} and $\Theta(b)$ in \eqref{eq:class-lowbd},
the power of the test $\psi$ in \eqref{eq:chen-test} satisfies
\[
\lim_{n\to\infty} \inf_{\Theta(b)} \mathsf{E}_{\Sigma}\psi= 1- \Phi
\biggl(z_{1-\alpha} - \frac{b^2}{2} \biggr) > \alpha.
\]
Moreover, for $b_n\to\infty$, $\lim_{n\to\infty} \inf_{\Theta
(b_n)} \mathsf{E}_{\Sigma}\psi= 1$.
\end{theorem}

Theorem \ref{thm:uppbd} shows that the test $\psi$ can distinguish
between the null \eqref{eq:H0-cov} and the alternative \eqref
{eq:H1-cov} with power tending to $1$ when $b=b_n\to\infty$.
Comparing with the lower bound given in Theorem~\ref{thm:lowbd}, the
test $\psi$ is rate-optimal when $p/n$ is bounded.
When $\llVert  \Sigma-I \rrVert_F \asymp\sqrt{p/n}$, the proof of Theorem
\ref{thm:uppbd} essentially shows that the power of $\psi$ is also
monotone increasing in $\llVert  \Sigma-I \rrVert_F$.

To prove Theorem \ref{thm:uppbd}, we first notice that the second
claim is a direct consequence of the first one. Indeed, if the first
claim is true, then
for any fixed constant $b>0$,
\[
\liminf_{n\to\infty} \inf_{\Theta(b_n)} \mathsf{E}_\Sigma\psi \geq
\lim_{n\to\infty} \inf_{\Theta(b)}\mathsf{E}_\Sigma\psi= 1- \Phi
\biggl(z_{1-\alpha} - \frac{b^2}{2} \biggr).
\]
Because the above inequality holds for any $b$, we obtain $\liminf_{n\to\infty}\inf_{\Theta(b_n)} \mathsf{E}_\Sigma\psi\geq1$. On the
other hand, $\psi\leq1$ and so $\limsup_{n\to\infty}\inf_{\Theta
(b_n)} \mathsf{E}_\Sigma\psi\leq1$. This leads to the second claim.

Turn to the proof of the first claim, we divide $\Theta(b)$ into two
disjoint subsets $\Theta(b) = \Theta(b,B)\cup\Theta(B)$, where
\begin{eqnarray}\label{eq:class}
 \Theta(b, B) &=& \bigl\{\Sigma\dvt b\sqrt{p/n}\leq
\llVert \Sigma-I \rrVert_F < B\sqrt{p/n} \bigr\} ,
\nonumber\\[-8pt]\\[-8pt]
\Theta(B) &= &\bigl\{\Sigma\dvt \llVert \Sigma-I \rrVert_F \geq B
\sqrt{p/n} \bigr\}. \nonumber
\end{eqnarray}
Here, $B$ is a sufficiently large constant, the choice of which depends
only on $\alpha$ and $b$, but not on $n$ or $p$.
We
employ different proof strategies on the two subsets. On $\Theta(B)$,
Chebyshev's inequality readily shows that
\[
\inf_{\Theta(B)}\mathsf{E}_{\Sigma}\psi> 1 -\Phi \biggl(z_{1-\alpha
}-
\frac{b^2}{2} \biggr).
\]
Turn to $\Theta(b,B)$. On this subset, Proposition \ref
{prop:T-n-MGCLT} then plays the key role in obtaining a uniform
approximation to the power function $\mathsf{E}_\Sigma\psi$ by the normal
distribution function $\Phi(z_{1-\alpha} - \frac{\llVert  \Sigma -I
\rrVert_F^2}{2p/n})$, which in turn leads to the final claim. For a detailed
proof, see Section~\ref{sec:pf-uppbd}.

\begin{remark}
(a). When $p/n$ is bounded, the conclusion of the theorem matches the
lower bound in Theorem \ref{thm:lowbd}. However, the result here holds
even when $p/n$ is unbounded.

(b). It can be seen from the proof of Theorem \ref{thm:uppbd} that the
simple expression
\[
\Phi \biggl(z_{1-\alpha} - \frac{\llVert  \Sigma-I \rrVert_F^2}{2p/n} \biggr)
\]
gives good approximation to the power of the test $\psi$ defined in
\eqref{eq:chen-test} at any $\Sigma$ of interest in practice,
because the approximation works well until the power of the test is
extremely close to $\alpha$ or $1$.
\end{remark}

%s3.4 ###
\subsection{Power comparison with the corrected LRT}
\label{sec:CLRT}

In the classical asymptotic regime where $p$ is fixed and $n\to\infty
$, the likelihood ratio test (LRT) is one of the most commonly used
tests. In the high-dimensional setting where both $n$ and $p$ are large
and $p<n$, Bai \textit{et al.} \cite{baji09} showed that the LRT is not well
behaved as the chi-squared limiting distribution under $H_0$ no longer holds.

For testing \eqref{eq:H0-cov}, when $p<n$ and $p/n\to c\in(0,1)$, Bai
\textit{et al.} \cite{baji09} proposed a corrected likelihood ratio test (CLRT)
with the test statistic $\mathit{CLR}_n$ given in \eqref{eq:clrt}. It was shown
that the test statistic $\mathit{CLR}_n \Rightarrow N(0,1)$ under $H_0$ and this
leads to an asymptotically level $\alpha$ test by rejecting $H_0$ when
$\mathit{CLR}_n > z_{1-\alpha}$. It was shown that the CLRT significantly
outperforms the LRT when both $n$ and $p$ are large and $p< n$.
Recently, Jiang \textit{et al.} \cite{jiji10} also considered the CLRT and
showed that the above limit holds even when $p/n\to1$.

It is interesting to compare the power of the CLRT with that of the
test defined in \eqref{eq:chen-test}. Note that the test given in
\eqref{eq:chen-test} is always well defined, but the CLRT is only
properly defined in the asymptotic regime where $p<n$ and $p/n\to c\in
(0, 1]$.
The following result shows that the power of the CLRT is uniformly
dominated by that of $\psi$ given in \eqref{eq:chen-test} over the
entire asymptotic regime under which the CLRT is applicable.

% \nb{is it possible to include $p/n\to1$?}

\begin{proposition}
\label{prop:CLRT-compare}
% Let $\Theta(b)$ be defined as in \eqref{eq:class-lowbd} with $0<b <
%1$.
Suppose that as $n\to\infty$, $p\to\infty$ with $p<n$ and $p/n\to
c\in(0,1]$. Let $\mathbb{C}(\tau) = \{\Sigma\dvt \llVert  \Sigma- I
\rrVert_F
= \tau\sqrt{p/n} \}$.
Then
for $\psi$ in \eqref{eq:chen-test} and the corrected LRT $\phi_{\mathrm{CLR}}$, we have
\[
\lim_{n\to\infty} \inf_{\mathbb{C}(\tau)} \mathsf{E}_{\Sigma
}\psi>
\limsup_{n\to\infty} \inf_{\mathbb{C}(\tau)} \mathsf{E}_{\Sigma
}
\phi_{\mathrm{CLR}}\qquad  \mbox{for all $\tau\in(0,1)$.}
\]
Moreover, for $\Theta(b)$ in \eqref{eq:class-lowbd} with $b\in(0,1)$,
\[
\lim_{n\to\infty} \inf_{\Theta(b)} \mathsf{E}_{\Sigma}\psi>
\limsup_{n\to\infty} \inf_{\Theta(b)} \mathsf{E}_{\Sigma}
\phi_{\mathrm{CLR}}.
\]
\end{proposition}

Hence, the CLRT is sub-optimal whenever it is properly defined. A proof
of Proposition \ref{prop:CLRT-compare} is given in Section \ref
{sec:proof-clrt-compare}.

%f1 ###
\begin{figure}

\includegraphics{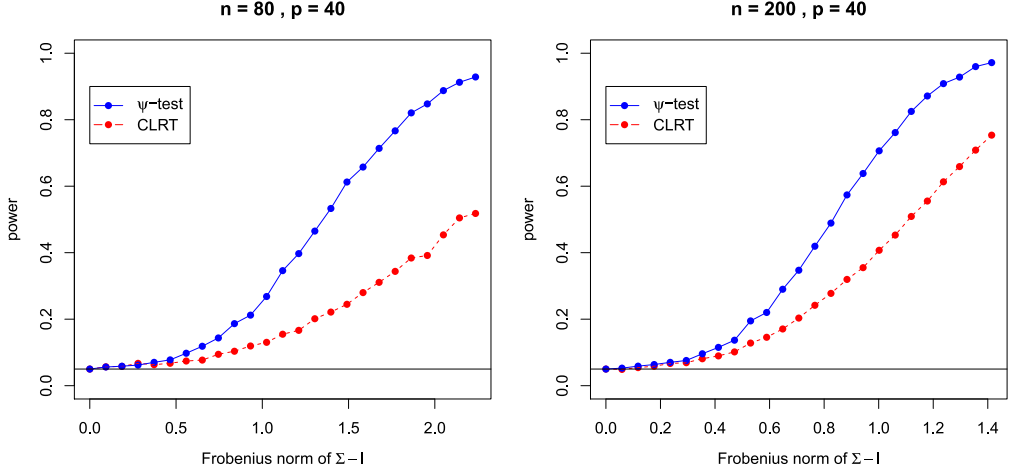}

\caption{Power curves of the $\psi$ test (solid) and the CLRT (dashed) under the equi-correlation
alternative. Each dot is obtained from $5000$ repetitions, and the
curves are then obtained via linear interpolation.}\label{fig:unif}
\end{figure}

%s4 ###
\section{Numerical experiments}
\label{sec:numeric}

In this section, a small simulation study is carried out to compare the
power of the test $\psi$ defined in \eqref{eq:chen-test} with that of
the CLRT under two specific alternatives.

The first alternative is the equi-correlation matrix $\Sigma= (\sigma_{ij})$, where for $\rho\in(0,1)$,
\[
\sigma_{ij} = %
\cases{ 1, &\quad  $i=j$,
\cr
\rho, &\quad  $i\neq j$. } %
\]
%
% \[
% \Sigma= \begin{bmatrix}
% 1 & \rho& \cdots& \rho\\
% \rho& \ddots& \ddots& \vdots\\
% \vdots& \ddots& \ddots& \rho\\
% \rho& \cdots& \rho& 1
% \end{bmatrix},  \mbox{with}  \rho> 0.
% \]
Figure \ref{fig:unif} shows how the power functions of the $\psi$
test and the CLRT grow with $\|\Sigma-I\|_F$ when $p = 40$ and $n =
80$ or $200$. For both tests, the significance levels are fixed at
$\alpha= 0.05$.
To make a fair comparison, the $95$th percentiles of the null
distributions of both test statistics are obtained via simulation
instead of using those of the asymptotic normal distributions.
From Figure \ref{fig:unif}, it is clear that the $\psi$ test is more
powerful than the CLRT for both $(n,p)$ configurations. The difference
between the powers is smaller when $n/p$ is larger. This is not
surprising, because the LRT is a powerful test in the ``large $n$,
small~$p$'' regime.

\begin{figure}

\includegraphics{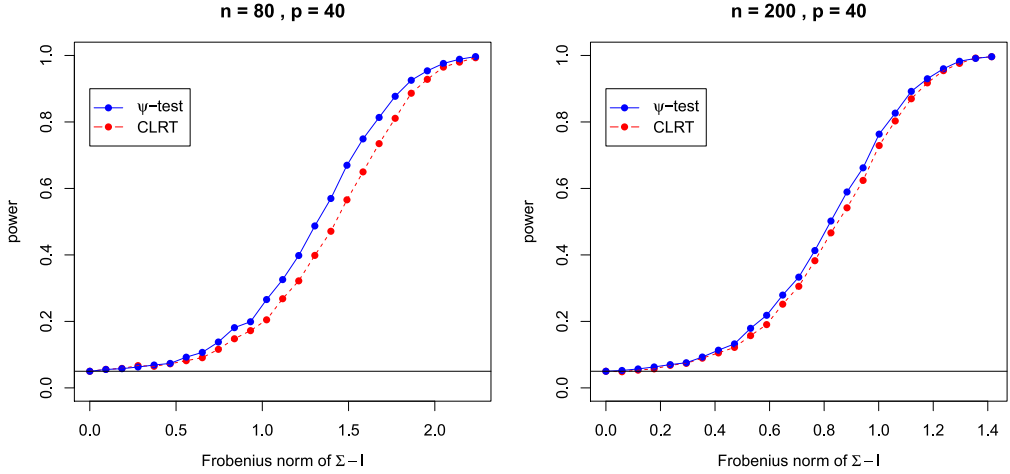}

\caption{Power curves of the $\psi$ test (solid) and the CLRT (dashed) against the tridiagonal alternative.
Each dot is obtained from $5000$ repetitions, and the curves are then
obtained via linear interpolation.}\label{fig:trid}
\end{figure}

The second alternative is the tridiagonal matrix $\Sigma= (\sigma_{ij})$, where for $\rho\in(0,1)$,
\[
\sigma_{ij} = %
\cases{ 1, &\quad  $i=j$,
\cr
\rho, &\quad  $|i-j|=1$,
\cr
0, & \quad $|i-j|>1$. } %
\]
Figure \ref{fig:trid} shows how the power functions of the $\psi$
test and the CLRT grow with $\|\Sigma-I\|_F$ for the tridiagonal
alternative. All the other setups remain unchanged. Here, the power of
the $\psi$ test still dominates, while the difference in power between
the two tests is smaller.

%f2 ###

%s5 ###
\section{Discussion}
\label{sec:disc}

% In this paper, we have considered the problem of minimax rate test
%for $H_0:\covmat= I$ in high dimensions where the distance metric is
%the Frobenius distance $\fnorm{\covmat-I}$. When the dimension to
%sample ratio $p/n$ is bounded, a lower bound of rate $O(\sqrt{p/n})$
%is established. Moreover, we show that there exists a test whose
%asymptotic power is uniformly greater than the significance level when
%the alternative is separated at rate $O(\sqrt{p/n})$ when both $n$ and
%$p$ tend to infinity, regardless of whether $p/n$ is bounded. Last but
%not least, the test uniformly dominates the corrected LRT in the
%minimax sense in all the high-dimensional settings where the latter is
%appropriate.

We have focused in the present paper on testing the hypotheses under
the Frobenius norm. The technical arguments developed in this paper can
also be used for testing under other matrix norms. Consider, for
example, testing \eqref{eq:H0-cov} against the following composite
alternative hypothesis
\[
% \label{eq:H1-cov-sp}
H_1\dvt \Sigma\in\Theta,\qquad  \mbox{where }  \Theta=
\Theta_n = \bigl\{\|\Sigma- I\|_s\geq\varepsilon_n
\bigr\}.
\]
Here $\|A\|_s$ is the spectral norm defined by $\|A\|_s = \max_{\|x\|
_2 =1} \|Ax\|_2$. Define
\begin{equation}
\label{eq:-sepctral-class-lowbd} \Theta_s(b) = \bigl\{\Sigma\dvt \|\Sigma-
I\|_s \geq b\sqrt{p/n} \bigr\}.
\end{equation}
Then the same lower bound holds for $\Theta_s(b)$. To be more precise,
we have the following result.
%That is, Theorem \ref{thm:lowbd} holds with $\Theta(b)$ replaced by $

\begin{theorem}
\label{thm:lowbd-spectral}
Let $0 < \alpha< \beta< 1$. Suppose that as $n\to\infty$, $p\to
\infty$ and $p/n\leq\kappa$ for some constant $\kappa<\infty$ and
all $n$. Then there exist a constant $b = b(\kappa,\beta-\alpha) <
1$, such that for any test $\phi$ with significance level $0<\alpha
<1$ for testing $H_0\dvt \Sigma= I$,
\[
\limsup_{n\to\infty}\inf_{\Sigma\in\Theta_s(b)} \mathsf {E}_{\Sigma
}\phi< \beta.
\]
\end{theorem}

The proof of Theorem \ref{thm:lowbd-spectral} is analogous to that of
Theorem \ref{thm:lowbd}. We believe that the rate of $\sqrt{p/n}$ in
the lower bound is sharp. It is however unclear which test is optimal
against the alternative \eqref{eq:-sepctral-class-lowbd} under the
spectral norm. Obtaining a matching upper bound for a practically
useful test is an interesting project for future research.

The results in the current paper also shed light on the problem of
testing for independence, that is, $H_0\dvt R = I$, where $R$ is the
population correlation matrix.
% Note that all the covariance matrices in the least favorable
%configuration $\Theta^*(b)$ \eqref{eq:least-favor-class} have diagonal
%elements all equal to $1$. Thus, they are also correlation matrices.
%Therefore,
Following Remark \ref{correlation.remark}, the proof of Theorems \ref
{thm:lowbd} and \ref{thm:lowbd-spectral} can be used directly to
establish the same lower bound results on testing the correlation matrix.
% \label{thm:lowbd-corr}
% Let $\tilde\Theta(b) = \{R: \fnorm{R - I}\geq b\sqrt{p/n} \}$, and $0
%< \alpha< \beta< 1$. Suppose that as $n\to\infty$, $p\to\infty$ and
%$p/n\leq\kappa$ for some constant $\kappa<\infty$ and all $n$. Then
%there exist a constant $b = b(\kappa,\beta-\alpha) < 1$, such that for
%any test $\phi$ with significance level $0<\alpha<1$ for testing $H_0:
%R = I$,
% \[
% \limsup_{n\to\infty}\inf_{\covmat\in\tilde\Theta(b)} \Ex_{\covmat}
% \]
% if we split the sample into two halves, and use one half to estimate
%the coordinate-wise variances, we can also construct a test statistic
%similar to $V_n$ \eqref{eq:V-n} which would be uniformly powerful over
%any alternative with $\fnorm{R-I}\geq b\sqrt{p/n}$. Since this test is
%only of theoretical interest, we omit the details here.
%Obtaining a matching upper bound for a practically useful test is an
%interesting project for future research.

% \nb{comment on Onatski, et al (2011)}

Onatski \textit{et al.} \cite{onmo11} also studied the hypothesis testing
problem \eqref{eq:H0-cov}, but their attention is restricted to
testing against alternatives that are rank one perturbations to the
identity matrix. That is, under the alternative $H_1$ the covariance
matrix belongs to the set $\Theta_h = \{I + h vv'\dvt \|v\|_2=1\}$. The
asymptotic regime is restricted to $p/n\to c\in(0,\infty)$. In this
asymptotic regime, Theorem 7 in \cite{onmo11} gives a lower bound
result analogous to Theorem \ref{thm:lowbd}. However, it does not
cover the case when $p/n\to0$, nor can it be extended to the case of
testing correlation matrices. In addition, we notice that though the
result in \cite{onmo11} enables one to study the asymptotic power of
all the eigenvalue-based tests on each $\Theta_h$ when $p/n\to c\in
(0,\infty)$, it does not give a minimax claim as we did in Theorem
\ref{thm:uppbd}.
%Nor does it apply to the cases where $p/n\to0$ or $\infty$.

The results in this paper also raised a number of interesting questions
for future research. One example is the testing of equality of two
covariance matrices based on the independent random samples $X_1,\dots
, X_{n_1}\stackrel{\mathrm{i.i.d.}}{\sim} N_p(\mu_1, \Sigma_1)$ and $Y_1,\dots
, Y_{n_2}\stackrel{\mathrm{i.i.d.}}{\sim} N_p(\mu_2, \Sigma_2)$. The validity
of many commonly used statistical procedures including the classical
Fisher's linear discriminant analysis requires the assumption of equal
covariance matrices. So it is of interest to test $H_0\dvt\Sigma_1 =
\Sigma_2$. Motivated by an unbiased estimator of the Frobenius norm
of $\Sigma_1-\Sigma_2$, Chen and Li \cite{chenli12} proposed a test
using a linear combination of $U$-statistics and studied its power. Cai
\textit{et al.} \cite{cai11} introduced a test based on the maximum of the
standardized differences between the entries of the two sample
covariance matrices. The test is shown to be powerful against sparse
alternatives and robust with respect to the population distributions.
However, the optimality of the two-sample tests has not been well
studied. This is an important topic for future research that is of both
theoretical and practical interest.

In the present paper, no structural assumption is imposed on the
alternative class of the covariance matrices such as sparsity or
bandedness. An optimal test against a structured alternative is
potentially very different from the test \eqref{eq:chen-test}
considered here. Recently, Cai and Jiang \cite{caji11} considered
testing the null hypothesis that $\Sigma$ is a banded matrix and
introduced a test based on the coherence of a random matrix.
Xiao and Wu \cite{xiwu11} proposed a test for testing $H_0\dvt \Sigma=
I$ against sparse alternatives. Their test is based on the
maximum of the standardized entries of the sample covariance matrix.
The limiting null distribution is shown to be a type I extreme value
distribution, the power of the test is not analyzed. It is interesting
to investigate the optimality of these testing problems with structured
alternatives.

%s6 ###
\section{Proofs}
\label{sec:proof}

In this section, we prove Theorems \ref{thm:lowbd}, \ref{thm:uppbd}
and Propositions \ref{prop:T-n-MGCLT} and \ref{prop:CLRT-compare}.

%s6.1 ###
\subsection{\texorpdfstring{Proof of Theorem \protect\ref{thm:lowbd}}{Proof of Theorem 1}}
\label{sec:pf-lowbd}

Recall that $P_0$ is the probability measure when $X_1,\dots
X_n\stackrel{\mathrm{i.i.d.}}{\sim}
N_p(0,I)$ and $P_v$ is the probability measure when $X_1,\dots,
X_n\stackrel{\mathrm{i.i.d.}}{\sim}N_p(0,\Sigma_v)$. In addition, $P_1 = \frac
{1}{2^p}\sum_{v\in\{\pm1\}^p} P_v$ is the average measure of the $P_v$'s. Let
$f_0$ and $f_1$ be the density functions of $P_0$ and $P_1$, respectively.
By the discussion following Theorem \ref{thm:lowbd}, we could prove
Theorem \ref{thm:lowbd} by showing that $\int f_1^2/f_0 - 1 \leq
4(\beta-\alpha)^2$.

After some basic calculation (see Appendix \ref{app:proof-lowbd} for
details), we obtain that if $b < b_0(\kappa)$ such that
\begin{equation}
\label{eq:b-cond} b<1\quad   \mbox{and}\quad   \frac{bp}{\sqrt{n(p-1)}} <
\frac{1}{\sqrt{2}},
\end{equation}
then
\begin{equation}
\label{eq:l2-rep} \int\frac{f_1^2}{f_0} = \frac{(1-a^2)^{n-np/2}}{ [1+(p-1)a^2 ]^{n}} \mathsf {E}
\biggl[1- \biggl(\frac{pa}{1+(p-1)a^2} \biggr)^2 \biggl(
\frac{\mathbf
{1}'V}{p} \biggr)^2 \biggr]^{-n/2}.
\end{equation}
Here, the expectation is taken w.r.t.~$V = (V_1,\dots, V_p)'$ where
the $V_j$'s are i.i.d.~Rademacher random variables which take values
$\pm1$ with equal probability.

Note that \eqref{eq:b-cond} and $p/n \leq\kappa$ implies
\[
\biggl(\frac{pa}{1+(p-1)a^2} \biggr)^2 \leq\frac{1}{2}.
\]
%
% \underline{If $p/n\in[0, y_U]$ as $n\to\infty$}, then for $n\geq
%n_0(y_U)$, the above choice of $a$ ensures that the $t<1/2$ in
Also note that $(\mathbf{1}'V/p)^2\in[0,1]$. Thus, let $\tilde
{b}_{np} =
(\frac{pa}{1+(p-1)a^2})^2$, and $(\mathbf{1}'V/p)^2 = \xi_p$, we have
\begin{eqnarray*}
\mathsf{E}(1-\tilde{b}_{np}\xi_p)^{-n/2} & =&
\mathsf{E} \bigl[(1-\tilde {b}_{np}\xi_p)^{-1/(\tilde{b}_{np}\xi_p)}
\bigr]^{n\tilde{b}_{np}\xi
_p/2}
\\
& \leq&\mathsf{E}\exp \biggl(\frac{\log{4}}{2} n\tilde{b}_{np}
\xi_p \biggr) \qquad \bigl(0\leq(1-x)^{1/x}\leq4,\mbox{ for all }
x\in[0,1/2]\bigr)
\\
& \leq&\mathsf{E}\exp \biggl(\frac{\log{4}}{2} \frac{b^2p^2}{p-1}
\xi_p \biggr)\qquad  \bigl(\tilde{b}_{np}\leq
p^2b^2/\bigl[n(p-1)\bigr]\bigr).
\end{eqnarray*}
For $\xi_p$, Hoeffding's inequality \cite{hoef63}, applied to
Rademacher variables, yields
\[
\mathsf{P} (\xi_p \geq\lambda )\leq2\mathrm{e}^{-2p\lambda}\qquad    \mbox
{for all }\lambda> 0.
\]
Thus, we obtain
\begin{eqnarray*}
\mathsf{E}\exp \biggl(\frac{\log{4}}{2} \frac{b^2p^2}{p-1} \xi_p
\biggr) & = &\int_0^\infty\mathsf{P} \biggl(\exp
\biggl(\frac{\log{4}}{2}  \frac
{b^2p^2}{p-1} \xi_p \biggr)\geq u
\biggr)\,\mathrm{d}u
\\
& = &1 + \int_1^\infty\mathsf{P} \biggl(
\xi_p\geq\frac{2\log
u}{\log
4}\frac{p-1}{b^2p^2} \biggr)\,\mathrm{d}u
\\
& \leq&1 + \int_1^\infty2\exp \biggl(-
\frac{4(p-1)\log{u}}{b^2p\log
4} \biggr)\,\mathrm{d}u
\\
& = &1 + \frac{2b^2 p \log{2}}{2(p-1)-b^2 p \log{2}}.
\end{eqnarray*}
Here, the last equality holds if $2(p-1)>b^2 p \log{2}$, which is
always true for large $p$ since $b < 1$.

% \begin{remark}
% The above calculation works for the case $p\asymp n^\gamma$, $\gamma
%> 1$.
% \end{remark}

In addition, with $b$ satisfying \eqref{eq:b-cond}, when $n\to\infty$,
\[
\bigl(1-a^2\bigr)^{n-np/2}   \to \mathrm{e}^{b^2},\qquad
\bigl[1+(p-1)a^2\bigr]^n\to \mathrm{e}^{b^2}.
\]
Therefore, for large enough $n\geq n_0(\kappa)$,
\[
% \|P_1 - P_0\|_1^2 \leq
\int\frac{f_1^2}{f_0} - 1 \leq\frac{8b^2 p \log{2}}{2(p-1)-b^2 p \log{2}},
\]
which, for sufficiently small $b\leq b_0(\kappa,\beta-\alpha)$, is
no larger than $4(\beta-\alpha)^2$. This completes the proof.

%s6.2 ###
\subsection{\texorpdfstring{Proof of Proposition \protect\ref{prop:T-n-MGCLT}}{Proof of Proposition 3}}
\label{sec:pf-MGCLT}

Following the outline of proof after Proposition \ref{prop:T-n-MGCLT},
for $E_1$ and $E_2$ defined in \eqref{eq:E1-E2}, we complete the proof
below by showing that $E_1/\sigma^4_n(\Sigma) = \mathrm{O}(1/n)$ and
$E_2/\sigma^4_n(\Sigma) = \mathrm{O}(\operatorname{tr}(\Sigma^4)/\operatorname{tr}^2(\Sigma^2))$.

To this end, we start with some preliminaries. Throughout the proof,
$\mathsf{E}$ and $\operatorname{\mathsf{Var}}$ are used as
abbreviations for $\mathsf{E}_{\Sigma}$ and
$\operatorname{\mathsf{Var}}_{\Sigma}$, respectively.
% $T_n$ and $\sigma_n^2$ are abbreviations of $T_n(\covmat)$ and $
Recall the martingale representation \eqref{eq:T_n-MG},
% of $T_n = \sum_{k=1}^n D_{nk}$,
where each martingale difference $D_{nk}$ has the explicit expression
\eqref{eq:T_n-MG-diff}.
For $D_{nk}$, its conditional variance is
\begin{eqnarray}
\label{eq:sigma_nk} %
\sigma_{nk}^2  &=&\mathsf{E}_{k-1}
\bigl[D_{nk}^2\bigr]\nonumber\\
 &=& \frac{8}{n^2(n-1)^2}
\operatorname{tr}(Q_{k-1}\Sigma Q_{k-1}\Sigma ) +
\frac
{16}{n^2(n-1)}\operatorname{tr}\bigl(Q_{k-1}\Sigma^3
\bigr)
\\
&&{}  - \frac{16}{n^2(n-1)}\operatorname{tr}\bigl(Q_{k-1}
\Sigma^2\bigr) + \frac{8}{n^2}\operatorname{tr} \bigl(
\Sigma^2(\Sigma- I)^2 \bigr). \nonumber
\end{eqnarray}
Detailed derivation of \eqref{eq:T_n-MG-diff} and \eqref{eq:sigma_nk}
can be found in Appendix \ref{app:prop-rate-detail}.
With \eqref{eq:sigma_nk}, it is not difficult to verify that
\[
\mathsf{E}\bigl[\sigma_{nk}^2\bigr] = \frac{8(k-1)}{n^2(n-1)^2}
\bigl[\operatorname{tr}^2\bigl(\Sigma^2\bigr) +
\operatorname{tr}\bigl(\Sigma^4\bigr) \bigr] + \frac{8}{n^2}
\operatorname {tr} \bigl(\Sigma^2(\Sigma- I)^2 \bigr),
\]
and that $\sigma_n^2 = \operatorname{\mathsf{Var}}(T_n) = \sum_{k=1}^n \mathsf{E}[\sigma_{nk}^2]$. Last but not least, we have for any $j>k$,
\begin{equation}
\label{eq:increasing-seq} \mathsf{E}_{k-1}\bigl[\sigma_{nj}^2
- \mathsf{E}\sigma_{nj}^2\bigr] = \sigma_{nk}^2
- \mathsf{E} \sigma_{nk}^2.
\end{equation}
Now, we turn to the studies of $E_1$ and $E_2$.

\textit{Term} $E_1$.
We begin with the first term.
% To understand it, we need the following identity: for any psd $B\in
% \begin{align}
% \label{eq:fourth-moment-id}
% \Ex\big[Z'BZ - \tr(B)\big]^4 = 48\tr(B^4) + 12\tr^2(B^2)\leq C
% \end{align}
Decompose the covariance matrix (as in \cite{chzh10}) as \mbox{$\Sigma=
\Gamma\Gamma'$}, with $\Gamma\in\mathbb{R}^{p\times p}$. Then, we
have the
representation
\begin{equation}
\label{eq:X-rep} X_i = \Gamma Z_i,\qquad
Z_i\stackrel{\mathrm{i.i.d.}} {\sim}N_p(0,I),\qquad   i=1,\dots, n.
\end{equation}
We further define
\[
% \label{eq:A}
A = \Gamma'\Gamma, \qquad  M_{k-1} =
\Gamma'\sum_{i=1}^{k-1}
\bigl(X_iX_i'-\Sigma\bigr)\Gamma= A\sum
_{i=1}^{k-1}\bigl(Z_iZ_i'-I
\bigr)A.
\]
With the above definition, \eqref{eq:T_n-MG-diff} can be rewritten as
\[
% \label{eq:D_nk}
D_{nk} = \frac{2}{n(n-1)} \bigl[Z_k'
M_{k-1} Z_k - \operatorname {tr}(M_{k-1}) \bigr] +
\frac{2}{n} \bigl[Z_k'\bigl(A^2-A
\bigr)Z_k - \operatorname{tr}\bigl(A^2-A\bigr) \bigr].
\]
Therefore, we obtain from the Cauchy--Schwarz inequality and Lemma \ref
{lemma:moments} that
\begin{eqnarray*}
\mathsf{E}\bigl[D_{nk}^4\bigr] & \leq&\frac{C}{n^4}
\mathsf{E} \bigl[Z_k'\bigl(A^2-A
\bigr)Z_k - \operatorname{tr} \bigl(A^2-A\bigr)
\bigr]^4 + \frac{C}{n^4(n-1)^4}\mathsf{E} \bigl[Z_k'M_{k-1}Z_k
- \operatorname{tr} (M_{k-1}) \bigr]^4
\\
& \leq&\frac{C}{n^4}\operatorname{tr}^2 \bigl(
\Sigma^2(\Sigma -I)^2 \bigr) + \frac
{C}{n^4(n-1)^4}
\mathsf{E} \bigl[\operatorname{tr}^2\bigl(M_{k-1}^2
\bigr) \bigr].
\end{eqnarray*}
For $\operatorname{tr}(M_{k-1}^2)$, we use the following lemma, the
proof of which
is given in Appendix \ref{app:prop-rate-detail}.
% \nb{double check the coefficient in front of $(k-1)(k-2)$!}

\begin{lemma}
\label{lemma:tr-M-k-1-square}
For $\operatorname{tr}(M_{k-1}^2)$, we have
\begin{eqnarray*}
\mathsf{E} \bigl[\operatorname{tr}\bigl(M_{k-1}^2\bigr)
\bigr] & =& (k-1) \bigl[\operatorname{tr}^2\bigl(\Sigma^2
\bigr)+\operatorname{tr} \bigl(\Sigma^4\bigr) \bigr],
\\
\operatorname{\mathsf{Var}} \bigl[\operatorname{tr}\bigl(M_{k-1}^2
\bigr) \bigr] & =& (k-1) \bigl[24\operatorname{tr}\bigl(\Sigma^8\bigr)
+ 16\operatorname{tr}\bigl(\Sigma^6\bigr)\operatorname{tr}\bigl(
\Sigma^2\bigr)+8\operatorname{tr}^2\bigl(
\Sigma^4\bigr) + 8\operatorname{tr}\bigl(\Sigma^4\bigr)
\operatorname{tr}^2\bigl(\Sigma^2\bigr) \bigr]
\\
&&{}  + 2(k-1) (k-2) \bigl[6\operatorname{tr}\bigl(\Sigma^8\bigr) +
2\operatorname{tr}^2\bigl(\Sigma^4\bigr) \bigr].
\end{eqnarray*}
\end{lemma}

For any sequences $\{a_n\}$ and $\{b_n\}$ of positive numbers, write
$a_n\lesssim b_n$ if $\limsup_{n\to\infty} a_n/\break b_n < \infty$.
Note that $\operatorname{tr}(\Sigma^{6})\leq\operatorname
{tr}(\Sigma^4)\operatorname{tr}(\Sigma^2)$ and
$\operatorname{tr}(\Sigma^{8})\leq\operatorname{tr}^2(\Sigma^4)$.
Since $\operatorname{tr}(\Sigma^4) = \mathrm{o}(\operatorname{tr}^2(\Sigma^2))$, Lemma~\ref
{lemma:tr-M-k-1-square} implies that
\[
\mathsf{E}\bigl[D_{nk}^4\bigr] \lesssim\frac{1}{n^{4}}
\operatorname{tr}^2 \bigl(\Sigma^2(\Sigma
-I)^2 \bigr) + \frac{k^2}{n^8}\operatorname{tr}^4
\bigl(\Sigma^2\bigr),
\]
and hence
\begin{equation}
\label{eq:E1-asymp} E_1 \lesssim\frac{1}{n^3}
\operatorname{tr}^2 \bigl(\Sigma^2(\Sigma
-I)^2 \bigr) + \frac{1}{n^5}\operatorname{tr}^4
\bigl(\Sigma^2\bigr).
\end{equation}

\textit{Term} $E_2$.
For $E_2$, we can simplify it as
\begin{eqnarray*}
E_2 & =& \mathsf{E} \Biggl[\sum_{k=1}^n
\bigl(\sigma_{nk}^2 - \mathsf {E}\sigma_{nk}^2
\bigr) \Biggr]^{2}
\\
& =& \mathsf{E} \Biggl[\sum_{k=1}^n \bigl(
\sigma_{nk}^2 - \mathsf {E}\sigma_{nk}^2
\bigr)^2 + 2\sum_{k=1}^{n-1}\sum
_{l=k+1}^n \bigl(\sigma_{nk}^2
- \mathsf{E}\sigma_{nk}^2 \bigr) \bigl(
\sigma_{nl}^2 - \mathsf {E}\sigma_{nl}^2
\bigr) \Biggr]
\\
& =& \sum_{k=1}^n \operatorname{
\mathsf{Var}}\bigl(\sigma_{nk}^2\bigr) + 2\sum
_{k=1}^{n-1}(n-k)\operatorname{\mathsf{Var}} \bigl(
\sigma_{nk}^2\bigr)
=\sum_{k=1}^n (2n-2k+1)\operatorname{
\mathsf{Var}}\bigl(\sigma_{nk}^2\bigr).
\end{eqnarray*}
Here, the second equality comes from \eqref{eq:increasing-seq}.

Note that $\operatorname{tr}(Q_{k-1}\Sigma Q_{k-1}\Sigma) =
\operatorname{tr}(M_{k-1}^2)$ and
$\operatorname{tr} (Q_{k-1}(\Sigma^3-\Sigma^2) ) =
\operatorname{tr}
(M_{k-1}(A^2-A) )$.
% \begin{align*}
% \tr(Q_{k-1}\covmat Q_{k-1}\covmat) = \tr(M_{k-1}^2),
% \tr\big(Q_{k-1}(\covmat^3-\covmat^2)\big) = \tr\big(M_{k-1}(A^2-A)
% \end{align*}
So, by \eqref{eq:sigma_nk}, there exist numeric constants $C$ and
$C'$, such that
\[
\operatorname{\mathsf{Var}}\bigl(\sigma_{nk}^2\bigr) \leq
\frac
{C}{n^4(n-1)^4}\operatorname{\mathsf{Var}} \bigl[\operatorname{tr}
\bigl(M_{k-1}^2\bigr) \bigr] + \frac{C'}{n^4(n-1)^2}
\operatorname{\mathsf {Var}} \bigl[\operatorname{tr} \bigl(M_{k-1}
\bigl(A^2-A\bigr) \bigr) \bigr].
\]
We have studied $\operatorname{\mathsf{Var}}[\operatorname
{tr}(M_{k-1}^2)]$ in Lemma \ref
{lemma:tr-M-k-1-square}. On the other hand, we have from Lemma \ref
{lemma:moments} that
\begin{eqnarray*}
\operatorname{\mathsf{Var}} \bigl[\operatorname {tr}\bigl(M_{k-1}
\bigl(A^2-A\bigr)\bigr) \bigr] & = &(k-1)\operatorname{\mathsf{Var}}
\bigl[\operatorname{tr} \bigl(AZZ'A\bigl(A^2-A\bigr)
\bigr) \bigr] = (k-1)\operatorname{\mathsf{Var}} \bigl[Z'
\bigl(A^4-A^3\bigr)Z \bigr]
\\
& =& (k-1) \bigl\{\mathsf{E} \bigl[ \bigl(Z'\bigl(A^4 -
A^3\bigr)Z \bigr)^2 \bigr] - \bigl[ \mathsf{E}
\bigl[Z'\bigl(A^4-A^3\bigr)Z\bigr]
\bigr]^2 \bigr\}
\\
& = &(k-1) \bigl[2\operatorname{tr} \bigl(\bigl(A^4 - A^3
\bigr)^2 \bigr) + \operatorname{tr}^2\bigl(A^4-A^3
\bigr) - \operatorname{tr}^2\bigl(A^4-A^3
\bigr) \bigr]
\\
& = &2(k-1)\operatorname{tr} \bigl(\Sigma^6(\Sigma-I)^2
\bigr)
\\
& \leq&2(k-1)\operatorname{tr}\bigl(\Sigma^4\bigr)\operatorname{tr}
\bigl(\Sigma^2(\Sigma-I)^2 \bigr).
\end{eqnarray*}
Since $\operatorname{tr}(\Sigma^{6})\leq\operatorname{tr}(\Sigma^4)\operatorname{tr}(\Sigma^2)$, $\operatorname{tr}
(\Sigma^{8})\leq\operatorname{tr}^2(\Sigma^4)$ and $\operatorname
{tr}(\Sigma^4) = \mathrm{o}(\operatorname{tr}^2(\Sigma^2))$, we obtain that
\[
\operatorname{\mathsf{Var}}\bigl(\sigma_{nk}^2\bigr)
\lesssim\frac
{k}{n^8}\operatorname{tr}\bigl(\Sigma^4\bigr)
\operatorname{tr}^2\bigl(\Sigma^2\bigr) +
\frac{k^2}{n^2}\operatorname{tr}^2\bigl(\Sigma^4\bigr)
+ \frac
{k}{n^6}\operatorname{tr} \bigl(\Sigma^4\bigr)
\operatorname{tr} \bigl(\Sigma^2(\Sigma-I)^2 \bigr),
\]
which leads to the bound
\begin{equation}
\label{eq:E2-asymp} E_2 \lesssim\frac{1}{n^3}\operatorname{tr}
\bigl(\Sigma^4\bigr)\operatorname {tr} \bigl(\Sigma^2(
\Sigma -I)^2 \bigr) + \frac{1}{n^4}\operatorname{tr}^2
\bigl(\Sigma^4\bigr)+\frac
{1}{n^5}\operatorname{tr} \bigl(
\Sigma^4\bigr)\operatorname{tr}^2\bigl(
\Sigma^2\bigr).
\end{equation}

\textit{Summing up}.
By \eqref{eq:var-T_n}, we have
\begin{equation}
\label{eq:sigma_n-asymp} \sigma_n^4 \asymp
\frac{1}{n^2}\operatorname{tr}^2 \bigl(\Sigma^2(
\Sigma -I)^2 \bigr) + \frac{1}{n^3}\operatorname{tr}^2
\bigl(\Sigma^2\bigr)\operatorname{tr} \bigl(\Sigma^2(
\Sigma-I)^2 \bigr) + \frac{1}{n^4}\operatorname{tr}^4
\bigl(\Sigma^2\bigr).
\end{equation}
Here and after, for any sequences $\{a_n\}$ and $\{b_n\}$ of positive
numbers, we write $a_n\asymp b_n$ if $a_n/b_n$ is bounded away from
both $0$ and $\infty$.
Thus, we obtain that
\[
\sigma_n^{-4}E_1 = \mathrm{O} \bigl(n^{-1}
\bigr),\qquad   \sigma_n^{-4}E_2 = \mathrm{O} \bigl(
\operatorname{tr}\bigl(\Sigma^4\bigr)/\operatorname{tr}^2
\bigl(\Sigma^2\bigr) \bigr). % + n^{-1}\tr^2(\covmat^2)/\tr^3(\covmat^4)\big).
\]
Plugging these estimates in Lemma \ref{lemma:MGCLT}, we complete the
proof. %\emph{uniform} convergence of
% \begin{align}
% \label{eq:b-e-bound}
% \sup_{x\in\RR}\big| \Pr(T_n\leq\sigma_n x) - \Phi(x) \big| \leq C
% \end{align}

%s6.3 ###
\subsection{\texorpdfstring{Proof of Theorem \protect\ref{thm:uppbd}}{Proof of Theorem 4}}
\label{sec:pf-uppbd}

Following the discussion after Theorem \ref{thm:uppbd}, we give below
the detailed proof of the first claim in the theorem. In particular, we
bound the power of the test separately on $\Theta(B)$ and $\Theta
(b,B)$, which are defined in \eqref{eq:class}.

% As a first step, we divide the set $\Theta(b)$ into two disjoint
%subsets $\Theta(b) = \Theta(b,B)\cup\Theta(B)$, where
% \begin{align}
% \label{eq:class}
% \Theta(b, B) = \Big\{\covmat:
% b\sqrt{p/n}\leq\fnorm{\covmat-I} < B\sqrt{p/n} \Big\},
% \Theta(B) = \Big\{\covmat:
% \fnorm{\covmat-I} \geq B\sqrt{p/n} \Big\}.
% \end{align}
% Here, $B$ is a large enough constant, the choice of which will be
%made explicit later. The most important property of $B$ is that it
%depends only on $\alpha$ and $b$, but not on $n$ or $p$.
%
% We
% employ different proof strategies on the two subsets. On $\Theta(B)$,
%we proceed heavy-handedly by using Chebyshev's inequality since the
%alternative class is far away enough from $H_0$. On $\Theta(b,B)$, we
%apply Proposition \ref{prop:T-n-MGCLT} to obtain uniform approximation
%to the power function by the normal distribution function.

\textit{Case 1:} $\Theta(B)$.
Here, we shall proceed heavy-handedly by using Chebyshev's inequality,
because the alternative class is sufficiently far away from $H_0$.

For any $\Sigma\in\Theta(B)$, there exists $\tau\geq B$,
s.t.~$\llVert  \Sigma-I \rrVert_F = \tau\sqrt{p/n}$. Suppose $B$
is large
enough s.t.~$\tau^2 \geq B^2 \geq3 z_{1-\alpha}$. Note that $\sigma_n(I) = (2p/n)(1+\mathrm{o}(1))$, and so
\[
\mathsf{E}_\Sigma T_n = \llVert \Sigma-I
\rrVert_F^2 = \frac{\tau^2}{2}\sigma_n(I)
\bigl(1+\mathrm{o}(1)\bigr) > z_{1-\alpha}  \sigma_n(I).
\]
Thus, we can use Chebyshev's inequality to bound the type II error of
$\psi$ at $\Sigma$ as the following:
\begin{eqnarray}\label{eq:type2bd}
1-\mathsf{E}_{\Sigma}\psi&=& \mathsf{P}_\Sigma
\bigl(T_n\leq z_{1-\alpha
}\sigma_n(I)\bigr) =
\mathsf{P}_{\Sigma}\bigl(T_n-\mathsf{E}_{\Sigma}T_n
\leq z_{1-\alpha}  \sigma_n(I)-\mathsf{E}_{\Sigma}T_n
\bigr)
\nonumber
\\
& \leq&\mathsf{P}_{\Sigma}\bigl(|T_n-\mathsf{E}_{\Sigma}T_n|
\geq |z_{1-\alpha
}\sigma_n(I) -\mathsf{E}_{\Sigma}T_n|
\bigr)
\\
& \leq&\frac{\operatorname{\mathsf{Var}}_{\Sigma
}(T_n)}{[z_{1-\alpha}\sigma_n(I) - \mathsf{E}_{\Sigma}T_n]^2}.\nonumber
\end{eqnarray}
For $\operatorname{\mathsf{Var}}_\Sigma(T_n) = \sigma_n^2(\Sigma
)$, we have its explicit
expression given in \eqref{eq:var-T_n}. Let $\lambda_{\max}(\Sigma
)$ denote the largest eigenvalue of $\Sigma$.
When $\llVert  \Sigma-I \rrVert_F = \tau\sqrt{p/n}$, we have
$\lambda_{\max
}(\Sigma)\leq1 + \tau\sqrt{p/n}$, and so
\[
\operatorname{tr}\bigl(\Sigma^2\bigr)  \leq p \biggl(1 +
\frac{\tau}{\sqrt {n}} \biggr)^2,\qquad   \operatorname{tr}\bigl(
\Sigma^2(\Sigma-I)^2\bigr) \leq\lambda_{\max
}^2(
\Sigma)\llVert \Sigma-I \rrVert_F^2 \leq
\frac{\tau^2 p}{n} \biggl(1+\tau\sqrt{\frac{p}{n}} \biggr)^2.
\]
Since $\operatorname{tr}(\Sigma^4)\leq\operatorname{tr}^2(\Sigma^2)$ and $\sigma_n^2(I) =
(4p^2/n^2)(1+\mathrm{o}(1))$, the above inequalities, together with \eqref
{eq:var-T_n}, lead to
\[
\frac{\sigma_n^2(\Sigma)}{\sigma_n^2(I)} \leq \biggl[2 \biggl(1 + \frac{\tau}{\sqrt{n}}
\biggr)^2 + \frac{2\tau^2}{p} \biggl(1 + \tau \sqrt{
\frac{p}{n}} \biggr)^2 \biggr] \bigl(1+\mathrm{o}(1) \bigr).
\]
Since $\tau^2/2-z_{1-\alpha}\geq\tau^2/6$ , there exists some
constant $C_\alpha$ depending only on $\alpha$, such that
\begin{eqnarray*}
\frac{\operatorname{\mathsf{Var}}_{\Sigma}(T_n)}{[z_{1-\alpha
}\sigma_n(I) - \mathsf{E}_{\Sigma}T_n]^2} & \leq&\frac{2 (1 + {\tau}/{\sqrt{n}})^2 +
({2\tau^2}/{p})(1 + \tau\sqrt{{p}/{n}})^2}{(\tau^2/2 - z_{1-\alpha
})^2}  \bigl(1+\mathrm{o}(1) \bigr)
\\
& \leq &C_\alpha \biggl[ \biggl(\frac{1}{\sqrt{\tau}} + \frac
{1}{\sqrt{n}}
\biggr)^4 + \frac{1}{\tau^2 p} + \frac{1}{n} \biggr].
\end{eqnarray*}

Note that all the $\mathrm{o}(1)$ terms in the above derivation are uniform over
$\Theta(B)$.
Therefore, given $\alpha$ and $b$, there exist a constant $B =
B(\alpha, b)$, such that
\begin{eqnarray}\label{eq:power-Theta-B}
\liminf_{n\to\infty} \inf_{\Theta(B)} \mathsf{E}_{\Sigma} \psi& \geq&
1 - \frac{\operatorname{\mathsf{Var}}_{\Sigma}(T_n)}{[z_{1-\alpha
}\sigma_n(I) - \mathsf{E}_{\Sigma}T_n]^2}
\nonumber
\\
& \geq&1 - C_\alpha \biggl[ \biggl(\frac{1}{\sqrt{B}} +
\frac
{1}{\sqrt{n}} \biggr)^4 + \frac{1}{B^2 p} + \frac{1}{n}
\biggr]
\\
& \geq&1- \Phi \biggl(z_{1-\alpha} - \frac{b^2}{2} \biggr) > \alpha.
\nonumber
\end{eqnarray}

\textit{Case 2:} $\Theta(b,B)$.
On this subset, we use Proposition \ref{prop:T-n-MGCLT} to obtain the
following uniform approximation to the power function by the normal
distribution function
\begin{equation}
\label{eq:power-approx} \sup_{\Theta(b,B)}\biggl\llvert \mathsf{E}_{\Sigma}
\psi- \Phi \biggl(z_{1-\alpha} - \frac{\llVert  \Sigma-I \rrVert_F^2}{2p/n} \biggr) \biggr\rrvert
\to0.
\end{equation}
If \eqref{eq:power-approx} is true, then we obtain
\[
\lim_{n\to\infty} \inf_{\Theta(b,B)} \mathsf{E}_{\Sigma}\psi =
\lim_{n\to\infty} \inf_{\Theta(b,B)}\Phi \biggl(z_{1-\alpha} -
\frac{\llVert  \Sigma-I \rrVert_F^2}{2p/n} \biggr) = 1- \Phi \biggl(z_{1-\alpha} - \frac{b^2}{2}
\biggr) > \alpha.
\]
Together with \eqref{eq:power-Theta-B}, this leads to the desired claim.

Turn to the proof of \eqref{eq:power-approx}. First, note that
uniformly on $\Theta(b,B)$, we have
\begin{equation}
\label{eq:covmat-local-behavior} p (1-B/\sqrt{n} )^2\leq
\operatorname{tr}\bigl(\Sigma^2\bigr)\leq p (1 + B/\sqrt{n}
)^2,
\end{equation}
and $\operatorname{tr}(\Sigma^4) \leq\lambda_{\max}^2(\Sigma)
\operatorname{tr}(\Sigma^2) \leq B^2 (p/n)  p  (1+B/\sqrt{n})^2$.
Therefore, as $n\to\infty$,
\begin{equation}
\label{eq:cond-verify} \sup_{\Theta(b,B)}|\operatorname{tr}\bigl(
\Sigma^2\bigr)/p-1|\to0,\qquad   \sup_{\Theta(b,B)}\operatorname{tr}
\bigl(\Sigma^4\bigr)/\operatorname {tr}^2\bigl(
\Sigma^2\bigr) \to0.
\end{equation}
So the condition of Proposition \ref{prop:T-n-MGCLT} is satisfied.
Next, we observe that
\begin{eqnarray*}
\mathsf{E}_{\Sigma}\psi& = &\mathsf{P}_{\Sigma}\bigl(T_n
> z_{1-\alpha} \sigma_n(I) \bigr) = \mathsf{P}_{\Sigma}
\biggl(\frac{T_n - \llVert  \Sigma-I \rrVert_F^2 }{\sigma_n(\Sigma)} \geq\frac{\sigma_n(I)}{\sigma_n(\Sigma)}z_{1-\alpha
} -
\frac{\llVert  \Sigma-I \rrVert_F^2}{\sigma_n(\Sigma)} \biggr)
\\
& =& \mathsf{P}_{\Sigma} \biggl(\frac{T_n - \mu_n(\Sigma) }{\sigma_n(\Sigma)} \geq
\frac{\sigma_n(I)}{\sigma_n(\Sigma)}z_{1-\alpha
} - \frac{\llVert  \Sigma-I \rrVert_F^2}{\sigma_n(\Sigma)} \biggr).
\end{eqnarray*}
%
% Here, the last equality holds as $T_n(\covmat) = T_n(I)-\fnorm{
Thus, Proposition \ref{prop:T-n-MGCLT} and \eqref{eq:cond-verify}
together imply that
\[
\sup_{\Theta(b,B)}\biggl\llvert \mathsf{E}_{\Sigma}\psi- \Phi \biggl(
\frac
{\sigma_n(I)}{\sigma_n(\Sigma)}z_{1-\alpha} - \frac{\llVert
\Sigma-I \rrVert_F^2}{\sigma_n(\Sigma)} \biggr) \biggr\rrvert
\to0.
\]
To complete the proof, what is left to be verified is that
\begin{equation}
\label{eq:sigma-rate} \sup_{\Theta(b,B)}\biggl\llvert \frac{\sigma_n^2(\Sigma)}{ 4p^2/n^2
}-1
\biggr\rrvert \to0,
\end{equation}
because it implies
$\sup_{\Theta(b,B)}\llvert \Phi (z_{1-\alpha} - \frac{\llVert
\Sigma-I \rrVert_F^2}{2p/n} ) - \Phi (\frac{\sigma_n(I)}{\sigma_n(\Sigma)}z_{1-\alpha} - \frac{\llVert  \Sigma-I \rrVert_F^2}
{\sigma_n(\Sigma)} ) \rrvert  \to0$, which together with the last
display before \eqref{eq:sigma-rate}, leads to \eqref{eq:power-approx}.
To show \eqref{eq:sigma-rate}, first recall the expression of $\sigma_n^2(\Sigma)$ in \eqref{eq:var-T_n}. By \eqref{eq:cond-verify}, we
obtain that the first term in \eqref{eq:var-T_n} is $4(p^2/n^2)
(1+\mathrm{o}(1))$ where $\mathrm{o}(1)$ is uniform on $\Theta(b,B)$. On the other hand,
\begin{eqnarray*}
\operatorname{tr}\bigl(\Sigma^2(\Sigma-I)^2\bigr) &\leq&
\lambda^2_{\max
}(\Sigma)\llVert \Sigma-I \rrVert_F^2
\\
&\leq& \biggl( 1+B\sqrt{\frac
{p}{n}} \biggr)^2
\frac
{Bp}{n} \leq C(B) \max \biggl(1,\frac{p}{n} \biggr)\cdot
\frac{p}{n}.
\end{eqnarray*}
Here, $C(B)$ is a constant depending only on $B$.
Therefore, we have that the second term in \eqref{eq:var-T_n} is of
order $\mathrm{o}(p^2/n^2)$ uniformly over $\Theta(b,B)$. Putting the two parts
together leads to \eqref{eq:sigma-rate}. This completes the proof.

%s6.4 ###
\subsection{\texorpdfstring{Proof of Proposition \protect\ref{prop:CLRT-compare}}{Proof of Proposition 5}}
\label{sec:proof-clrt-compare}

Fix any $\tau\in(0,1)$. At each dimension $p$, consider a single
point in $\mathbb{C}(\tau)$:
\[
\Sigma^* = I_{p\times p} + \tau\sqrt{\frac{p}{n}}  uu',
\]
where $u$ is an arbitrarily fixed unit vector in $\mathbb{R}^p$. Since
$\tau
< 1$, Proposition 10 in \cite{onmo11} leads to
\[
\lim_{n\to\infty} \mathsf{E}_{\Sigma^*} \phi_{\mathrm{CLR}} = 1- \Phi
\bigl(z_{1-\alpha} - h(\tau,c) \bigr)\qquad   \mbox{for }   h(\tau,c) =
\frac{\tau\sqrt{c} - \log(1+\tau\sqrt{c})}{\sqrt{-2\log(1-c)
- 2c}}.
\]
Note that for all $\tau>0$ and $c\in(0,1)$, $\tau^2/2 > h(\tau,c) >
0$. Therefore,
\begin{eqnarray*}
\lim_{n\to\infty} \inf_{\mathbb{C}(\tau)} \mathsf{E}_{\Sigma
}\psi& =& 1- \Phi
\biggl(z_{1-\alpha} - \frac{\tau^2}{2} \biggr)
\\
&>& 1- \Phi \bigl(z_{1-\alpha} - h(\tau,c) \bigr)
\\
&=& \lim_{n\to\infty} \mathsf{E}_{\Sigma^*} \phi_{\mathrm{CLR}}
\\
&\geq&\limsup_{n\to\infty} \inf_{\mathbb{C}(\tau)} \mathsf {E}_{\Sigma
}
\phi_{\mathrm{CLR}}.
\end{eqnarray*}
The proof of the second claim is obtained by replacing $\mathbb
{C}(\tau)$ with $\Theta(b)$ and $\tau$ with $b$ in the above arguments.

\begin{appendix}\label{app}

\section*{Appendix: Technical details}
\label{sec:tech}

\renewcommand{\thesubsection}{A.\arabic{subsection}}
\setcounter{subsection}{0}

%s6.1 ###
\subsection{\texorpdfstring{Proof details for Theorem \protect\ref{thm:lowbd}}{Proof details for Theorem 1}}
\label{app:proof-lowbd}

\renewcommand{\theequation}{\arabic{equation}}
\setcounter{equation}{30}

Here we give the calculation leading to \eqref{eq:l2-rep} in the proof
of Theorem \ref{thm:lowbd}.

Consider $\Theta^*(b)$ in \eqref{eq:least-favor-class}.
% \begin{align}
% \label{eq:least-favor-class}
% \Theta^*(b) = \left\{\covmat_{v} = (1-\frac{b}{\sqrt{n(p-1)}})I_{p
% \end{align}
For any $v\in\{\pm1\}^p$, we have $\llVert  \Sigma_v - I \rrVert_F
= b\sqrt {p/n}$, $\operatorname{diag}(\Sigma_v) = (1,\dots, 1)$, and for $a =
b/\sqrt{n(p-1)}$,
\begin{eqnarray}
\label{eq:covmat_v} %
\Sigma^{-1}_v & =& \frac{1}{1-a}I_{p\times p} -
\biggl[\frac{1}{1-a} - \frac{1}{1+(p-1)a} \biggr]\frac{1}{p}vv',
\nonumber\\[-8pt]\\[-8pt]
\det{\Sigma_v} & =& (1-a)^{p-1}\bigl[1+(p-1)a\bigr].\nonumber
\end{eqnarray}
Therefore, we have
\begin{eqnarray*}
f_0(x_1,\dots, x_n) & =& \frac{1}{(2\uppi)^{np/2}}
\exp \Biggl\{-\frac
{1}{2}\sum_{i=1}^n
x'_i x_i \Biggr\},
\\
f_1(x_1,\dots, x_n) & =& \frac{1}{2^p}
\sum_{v}\frac{1}{(2\uppi
)^{np/2} (\det{\Sigma_{v}})^{n/2}} \exp \Biggl\{-
\frac{1}{2}\sum_{i=1}^n
x'_i \Sigma_{v}^{-1}
x_i \Biggr\}
\\
& =& \frac{1}{(2\uppi)^{np/2}} \exp \Biggl\{-\frac{1}{2(1-a)} \sum
_{i=1}^n x'_i
x_i \Biggr\} \frac{1}{(1-a)^{n(p-1)/2} [1+(p-1)a]^{n/2}}
\\
&&{}  \times\frac{1}{2^p}\sum_{v} \exp
\Biggl\{ \frac
{1}{2p} \biggl[\frac{1}{1-a} - \frac{1}{1+(p-1)a} \biggr]
\sum_{i=1}^n \bigl(v'x_i
\bigr)^2 \Biggr\}.
\end{eqnarray*}
And so
\begin{eqnarray*}
\frac{f_1^2}{f_0} & =& \frac{1}{(1-a)^{n(p-1)} [1+a(p-1)
]^{n}} \times\frac{1}{(2\uppi)^{np/2}} \exp \Biggl
\{-\frac{1}{2} \biggl(\frac{1+a}{1-a} \biggr)\sum
_{i=1}^n x'_i
x_i \Biggr\}
\\
& &{} \times\frac{1}{2^{2p}} \Biggl\{ \sum_v
\exp \Biggl[\frac
{1}{p} \biggl(\frac{1}{1-a}-\frac{1}{1+(p-1)a}
\biggr)\sum_{i=1}^n \bigl(v'
x_i\bigr)^2 \Biggr]
\\
&&\hspace*{35pt}{}   + \sum_{v\neq u} \exp \Biggl[
\frac{1}{2p} \biggl(\frac
{1}{1-a}-\frac{1}{1+(p-1)a} \biggr) \Biggl(
\sum_{{i=1}}^n\bigl(v'x_i
\bigr)^2 + \sum_1^n
\bigl(u' x_i\bigr)^2 \Biggr) \Biggr] \Biggr
\}.
\end{eqnarray*}
Now we compute the integral. Fix any $v$, we have
\begin{eqnarray*}
& &\int\frac{1}{(2\uppi)^{np/2}}\exp \Biggl\{-\frac{1}{2} \biggl(
\frac
{1+a}{1-a} \biggr)\sum_{i=1}^n
x'_i x_i \Biggr\}\exp \Biggl\{
\frac
{1}{p} \biggl(\frac{1}{1-a} - \frac{1}{1+(p-1)a} \biggr)\sum
_{i=1}^n \bigl(v'x_i
\bigr)^2 \Biggr\}\,\mathrm{d}x
\\
&&\quad  = \biggl[\int\frac{1}{(2\uppi)^{p/2}}\exp \biggl\{-\frac{1}{2} \biggl(
\frac{1+a}{1-a} \biggr)x' x \biggr\}\exp \biggl\{
\frac{1}{p} \biggl(\frac{1}{1-a} - \frac{1}{1+(p-1)a} \biggr)
\bigl(v'x\bigr)^2 \biggr\}\,\mathrm{d}x \biggr]^n
\\
&&\quad  = \biggl(\frac{1+a}{1-a} \biggr)^{-np/2} \bigl[\mathsf{E}\exp (tY)
\bigr]^n
\\
&&\quad  = \biggl(\frac{1+a}{1-a} \biggr)^{-np/2}(1-2t)^{-n/2}\qquad
\mbox {(for $t\leq1/2$)},
\end{eqnarray*}
where $Y\sim\chi^2_{(1)}$, and
\begin{equation}
\label{eq:t} t = \biggl(\frac{1}{1-a} - \frac{1}{1+(p-1)a} \biggr) \biggl(
\frac
{1+a}{1-a} \biggr)^{-1}.
\end{equation}
In addition, fix any $v \neq u$, we have
\begin{eqnarray*}
&& \int\frac{1}{(2\uppi)^{np/2}}\exp \Biggl\{-\frac{1}{2} \biggl(
\frac
{1+a}{1-a} \biggr)\sum_{i=1}^n
x'_i x_i \Biggr\}\\
&&\qquad {}\times\exp \Biggl\{
\frac
{1}{2p} \biggl(\frac{1}{1-a} - \frac{1}{1+(p-1)a} \biggr) \Biggl(
\sum_{i=1}^n \bigl(v'x_i
\bigr)^2 + \bigl(u' x_i\bigr)^2
\Biggr) \Biggr\}\,\mathrm{d}x
\\
&&\quad =  \biggl[\int\frac{1}{(2\uppi)^{p/2}}\exp \biggl\{-\frac{1}{2} \biggl(
\frac{1+a}{1-a} \biggr)x' x \biggr\}\\
&&\hspace*{35pt} {}\times\exp \biggl\{
\frac{1}{2p} \biggl(\frac{1}{1-a} - \frac{1}{1+(p-1)a} \biggr) \bigl[
\bigl(v'x\bigr)^2 + \bigl(u' x
\bigr)^2 \bigr] \biggr\}\,\mathrm{d}x \biggr]^n.
\end{eqnarray*}
Let $X\sim N_p(0,I)$, $Z_1 =v' X/\sqrt{p}$, and $Z_2 = u' X/\sqrt {p}$. Then
\[
\left[\matrix{ Z_1
\cr
Z_2 } \right]%
\sim N_2 %
\biggl(\matrix{ %
\left[\matrix{ 0
\cr
0 } \right]%
, %
\left[\matrix{1 &
v'u/p
\cr
v'u/p & 1 } \right]%
}\biggr) %
,
\]
and so $Z_1^2 + Z_2^2 \stackrel{d}{=} (1+v'u/p)Y_1 + (1-v'u/p)Y_2$,
with $Y_i\stackrel{\mathrm{i.i.d.}}{\sim}\chi^2_{(1)}$. Therefore, the second
last display equals
\begin{eqnarray*}
&&\biggl(\frac{1+a}{1-a} \biggr)^{-np/2} \bigl[\mathsf{E}
\exp(t_1 Y_1 + t_2 Y_2)
\bigr]^n\\
&&\quad  = \biggl(\frac{1+a}{1-a} \biggr)^{-np/2}
(1-2t_1)^{-n/2}(1-2t_2)^{-n/2},
\end{eqnarray*}
where
\[
t_1 = \frac{t}{2} \biggl(1 + \frac{1}{p}v'u
\biggr),\qquad   t_2 = \frac{t}{2} \biggl(1 -
\frac{1}{p}v'u \biggr).
\]

Collecting terms, we obtain after some linear algebra that %\nb{to be
%double-checked!}
%
\begin{eqnarray*}
\int\frac{f_1^2}{f_0} & =& \frac{(1-a^2)^{n-np/2}}{ [1+(p-1)a^2 ]^{n}} \mathsf {E}_{V,U}
\biggl[1- \biggl(\frac{pa}{1+(p-1)a^2} \biggr)^2 \biggl(
\frac{V'U}{p} \biggr)^2 \biggr]^{-n/2}
\\
& =& \frac{(1-a^2)^{n-np/2}}{ [1+(p-1)a^2 ]^{n}} \mathsf{E}_{V} \biggl[1- \biggl(
\frac{pa}{1+(p-1)a^2} \biggr)^2 \biggl(\frac
{\mathbf{1}'V}{p}
\biggr)^2 \biggr]^{-n/2}.
\end{eqnarray*}
Here, both $V$ and $U$ have i.i.d.~Rademacher entries which take values
$\pm1$ with equal probability, and in the first expectation, $V$ and
$U$ are independent.

%s6.2 ###
\subsection{Variance of $T_n$}
%Proof of \eqref{eq:var-T_n}}
\label{app:var-T-n}

In this part, we establish the variance of $T_n$ given in \eqref{eq:var-T_n}.
We begin with a technical lemma, which is closely connected to \cite{chzh10}, Proposition A.1.

\begin{lemma}
\label{lemma:moments}
Let $Z_1,Z_2\stackrel{\mathrm{i.i.d.}}{\sim}N_p(0,I)$, and $M$, $N$ be two
$p\times p$ p.s.d.
matrices, then
\begin{eqnarray}\label{eq:cross-2nd-moment}
 \mathsf{E}\bigl[\bigl(Z_1'M Z_1\bigr)
\bigl(Z_1'N Z_1\bigr)\bigr] &=&
\operatorname {tr}(M)\operatorname{tr}(N) + 2\operatorname{tr}(MN);
\\\label {eq:4th-moment}
 \mathsf{E}\bigl[\bigl(Z_1'M Z_2
\bigr)^4\bigr] &=& 3\operatorname{tr}^2\bigl(M^2
\bigr) + 6\operatorname{tr}\bigl(M^4\bigr);
\\\label {eq:4th-moment-central}
 \mathsf{E}\bigl[\bigl(Z_1'M Z_1 -
\operatorname{tr}(M)\bigr)^4\bigr] &=& 48\operatorname {tr}
\bigl(M^4\bigr) + 12\operatorname{tr}^2
\bigl(M^2\bigr).
\end{eqnarray}
\end{lemma}

\begin{pf}
Denote the ordered eigenvalues of $M$ by $\lambda_1\geq\cdots\geq
\lambda_p$, and those of $N$ by $\mu_1\geq\cdots\geq\mu_p$. Let
$U_j\stackrel{\mathrm{i.i.d.}}{\sim}N(0,1)$, $j=1,\dots,p$. For \eqref
{eq:cross-2nd-moment}, we have
\begin{eqnarray*}
\mathsf{E}\bigl[\bigl(Z_1'M Z_1\bigr)
\bigl(Z_1'N Z_1\bigr)\bigr] & =& \mathsf{E}
\Biggl[ \Biggl(\sum_{j=1}^p
\lambda_j U_j^2 \Biggr) \Biggl(\sum
_{j=1}^p\mu_j U_j^2
\Biggr) \Biggr]
\\
& =& \sum_{j=1}^p\lambda_j
\mu_j\mathsf{E}\bigl[U_j^4\bigr]+\sum
_{j\neq l}^p \lambda_j
\mu_l\mathsf{E}\bigl[U_j^2\bigr]\mathsf{E}
\bigl[U_l^2\bigr]
\\
& =& 3\sum_{j=1}^p \lambda_j
\mu_j + \sum_{j\neq l}^p
\lambda_j\mu_l = \operatorname{tr}(M)
\operatorname{tr}(N) + 2\operatorname{tr}(MN).
\end{eqnarray*}

For \eqref{eq:4th-moment}, we define $V_j\stackrel{\mathrm{i.i.d.}}{\sim
}N(0,1)$, $j=1,\dots,p$,
which are independent from the $U_j$'s. Then
\begin{eqnarray*}
\mathsf{E}\bigl[\bigl(Z_1'M Z_2
\bigr)^4\bigr] & =& \mathsf{E} \Biggl[ \Biggl(\sum
_{j=1}^p \lambda_j U_j
V_j \Biggr)^4 \Biggr]
\\
& = &\sum_{j=1}^p\lambda_j^4
\mathsf{E}\bigl[U_j^4\bigr]\mathsf{E}
\bigl[V_j^4\bigr] + {4\choose2}\sum
_{j\neq l}^p\lambda_j^2
\lambda_l^2 \mathsf{E}\bigl[U_j^2
\bigr]\mathsf {E}\bigl[U_l^2\bigr]\mathsf{E}
\bigl[V_j^2\bigr]\mathsf{E}\bigl[V_l^2
\bigr]
\\
& =& 9\sum_{j=1}^p\lambda_j^4
+ 6\sum_{j\neq l}^p\lambda_j^2
\lambda_l^2 = 3\operatorname{tr}^2
\bigl(M^2\bigr) + 6\operatorname{tr}\bigl(M^4\bigr).
\end{eqnarray*}

Finally, for \eqref{eq:4th-moment-central}, we have
\begin{eqnarray*}
\mathsf{E}\bigl[\bigl(Z_1'MZ_1 -
\operatorname{tr}(M)\bigr)^4\bigr] & =& \mathsf{E} \Biggl[ \Biggl(\sum
_{j=1}^p\lambda_j
\bigl(U_j^2-1\bigr) \Biggr)^4 \Biggr]
\\
& =& \sum_{j=1}^p \lambda_j^4
\mathsf{E}\bigl[\bigl(U_i^2-1\bigr)^4\bigr]
+ 6\sum_{j\neq
l}^p\lambda_j^2
\lambda_l^2 \mathsf{E}\bigl[\bigl(U_j^2-1
\bigr)^2\bigr] \mathsf {E}\bigl[\bigl(U_l^2-1
\bigr)^2\bigr]
\\
& =& 60\sum_{j=1}^p\lambda_i^4
+ 24\sum_{j\neq l}^p\lambda_j^2
\lambda_l^2 = 48\operatorname{tr}\bigl(M^4
\bigr) + 12\operatorname{tr}^2\bigl(M^2\bigr).
\end{eqnarray*}
This completes the proof of the lemma.
\end{pf}

In order to understand the variance of $T_n$, we need the following lemma.

\begin{lemma}
\label{lemma:var-cov-h}
For $X_1, X_2, X_3\stackrel{\mathrm{i.i.d.}}{\sim}N_p(0,\Sigma)$, we have
\begin{eqnarray*}
\operatorname{\mathsf{Var}} \bigl(h(X_1,X_2) \bigr) & =&
2 \bigl[\operatorname{tr}^2\bigl(\Sigma^2\bigr) +
\operatorname{tr}\bigl(\Sigma^4\bigr)\bigr] + 4 \operatorname{tr}
\bigl(\Sigma^2(\Sigma-I)^2 \bigr),
\\
\operatorname{\mathsf{Cov}} \bigl(h(X_1,X_2),
h(X_1,X_3) \bigr) & =& 2\operatorname{tr} \bigl(
\Sigma^2(\Sigma- I)^2 \bigr).
\end{eqnarray*}
\end{lemma}

\begin{pf}%[Proof of Lemma \ref{lemma:var-cov-h}]
For the variance, we first decompose it as
\[
\operatorname{\mathsf{Var}} \bigl(h(X_1,X_2) \bigr) =
\operatorname {\mathsf{Var}}\bigl(X_1'X_2
\bigr)^2 + 2\operatorname{\mathsf{Var}}\bigl(X_1'X_1
\bigr) - 4\operatorname{\mathsf{Cov}} \bigl(\bigl(X_1'X_2
\bigr)^2, \bigl(X_1'X_1\bigr)
\bigr).
\]
For $\operatorname{\mathsf{Var}}(X_1'X_2)^2 = \mathsf
{E}[(X_1'X_2)^4] - [\mathsf{E}(X_1'X_2)^2]^2$, we
have from \eqref{eq:4th-moment} that
\[
\mathsf{E}\bigl[\bigl(X_1'X_2
\bigr)^4\bigr] = \mathsf{E}\bigl[\bigl(Z_1'AZ_2
\bigr)^4\bigr] = 3\operatorname {tr}^2
\bigl(A^2\bigr) + 6\operatorname{tr}\bigl(A^4\bigr) = 3
\operatorname{tr}^2\bigl(\Sigma^2\bigr) + 6
\operatorname{tr}\bigl(\Sigma^4\bigr).
\]
On the other hand, we have
\[
\mathsf{E}\bigl[\bigl(X_1'X_2
\bigr)^2\bigr] = \mathsf{E}\bigl[Z_1' A
Z_2 Z_2' A Z_1\bigr] = \mathsf
{E}\bigl[\operatorname{tr}\bigl(AZ_2 Z_2' A
\bigr)\bigr] = \mathsf{E}\bigl[Z_2'A^2
Z_2\bigr] = \operatorname{tr}\bigl(A^2\bigr) =
\operatorname{tr}\bigl(\Sigma^2\bigr).
\]
Thus, we obtain $\operatorname{\mathsf{Var}}(X_1'X_2)^2 =
2\operatorname{tr}^2(\Sigma^2) + 6\operatorname{tr}(\Sigma^4)$. Similar type of calculation yields that
\[
\operatorname{\mathsf{Var}}\bigl(X_1'X_1
\bigr)  = 2\operatorname{tr} {\bigl(\Sigma^2\bigr)},\qquad
\operatorname{\mathsf{Cov}} \bigl(\bigl(X_1'X_2
\bigr)^2, \bigl(X_1'X_1\bigr)
\bigr) = 2\operatorname{tr} {\bigl(\Sigma^3\bigr)}.
\]
Assembling the pieces, we prove the variance formula.

For the covariance formula, the basic quantity to compute is
\begin{eqnarray*}
&& \mathsf{E}\bigl[\bigl(X_1'X_2
\bigr)^2 - \bigl(X_1'X_1\bigr)
- \bigl(X_2'X_2\bigr)\bigr] \bigl[
\bigl(X_1'X_3\bigr)^2 -
\bigl(X_1'X_1\bigr) -
\bigl(X_3'X_3\bigr)\bigr]
\\
&&\quad =   \mathsf{E}\bigl(X_1'X_2
\bigr)^2\bigl(X_1'X_3
\bigr)^2 - \mathsf {E}\bigl(X_1'X_1
\bigr) \bigl(X_1'X_3\bigr)^2 -
\mathsf{E} \bigl(X_2'X_2\bigr)\mathsf{E}
\bigl(X_1'X_3\bigr)^2
\\
&&\qquad {} - \mathsf{E}\bigl(X_1'X_2
\bigr)^2\bigl(X_1'X_1\bigr) +
\mathsf{E}\bigl(X_1'X_1
\bigr)^2 + \mathsf {E}\bigl(X_1'X_1
\bigr)\mathsf{E} \bigl(X_2'X_2\bigr)
\\
&&\qquad {} - \mathsf{E}\bigl(X_3'X_3\bigr)
\mathsf{E}\bigl[\bigl(X_1'X_2
\bigr)^2 - \bigl(X_1'X_1\bigr)
- \bigl(X_2'X_2\bigr)\bigr]
\\
&&\quad =   \mathsf{E}\bigl(X_1'X_2
\bigr)^2\bigl(X_1'X_3
\bigr)^2 - 2 \mathsf {E}\bigl(X_1'X_1
\bigr) \bigl(X_1'X_3\bigr)^2\\
&&\qquad {} -
2\mathsf{E} \bigl(X_1'X_1\bigr)\mathsf{E}
\bigl(X_1'X_2\bigr)^2 +
\mathsf{E}\bigl(X_1'X_1
\bigr)^2 + 3\bigl[\mathsf {E}\bigl(X_1'X_1
\bigr)\bigr]^2.
\end{eqnarray*}
First, we compute $\mathsf{E}(X_1'X_2)^2(X_1'X_3)^2$, for which we have
\begin{eqnarray*}
\mathsf{E}\bigl(X_1'X_2
\bigr)^2\bigl(X_1'X_3
\bigr)^2 & = &\mathsf{E}\bigl[\bigl(Z_1' A
Z_2\bigr)^2 \bigl(Z_1' A
Z_3\bigr)^2\bigr]
\\
& =& \mathsf{E}\bigl[\mathsf{E}\bigl[Z_2'AZ_1
Z_1' A Z_2|Z_1\bigr]
\mathsf {E}\bigl[Z_3'AZ_1
Z_1' A Z_3|Z_1\bigr]\bigr]
\\
& = &\mathsf{E}\bigl[\operatorname{tr}^2\bigl(AZ_1
Z_1'A\bigr)\bigr]
\\
& =& \mathsf{E}\bigl[\bigl(Z_1' A^2
Z_1\bigr)^2\bigr]
\\
& =& \operatorname{tr}^2\bigl(A^2\bigr) + 2
\operatorname{tr}\bigl(A^4\bigr) = \operatorname{tr}^2
\bigl(\Sigma^2\bigr) + 2\operatorname{tr}\bigl(\Sigma^4
\bigr).
\end{eqnarray*}
Next, we compute $\mathsf{E}(X_1'X_1)(X_1'X_3)^2$, for which we have
\begin{eqnarray*}
\mathsf{E}\bigl(X_1'X_1\bigr)
\bigl(X_1'X_3\bigr)^2 & =&
\mathsf{E}\bigl[\bigl(Z_1' A Z_1\bigr)
\bigl(Z_1' A Z_3\bigr)^2\bigr]
\\
& =& \mathsf{E}\bigl[\bigl(Z_1'A Z_1
\bigr)\mathsf{E}\bigl[ Z_3' A Z_1
Z_1' A Z_3 |Z_1\bigr]\bigr]
\\
& =& \mathsf{E}\bigl[\bigl(Z_1' A Z_1
\bigr) \bigl(Z_1' A^2 Z_1
\bigr)\bigr]
\\
& =& 2\operatorname{tr}\bigl(A^3\bigr) + \operatorname{tr}
\bigl(A^2\bigr)\operatorname {tr}(A) = 2\operatorname{tr}\bigl(
\Sigma^3\bigr) + \operatorname{tr}\bigl(\Sigma^2\bigr)
\operatorname{tr}(\Sigma).
\end{eqnarray*}
We further note that $\mathsf{E}(X_1'X_1)^2 = \mathsf{E}(Z_1' A
Z_1)^2 = 2\operatorname{tr}
(\Sigma^2)+\operatorname{tr}^2(\Sigma)$, that $\mathsf
{E}(X_1'X_2)^2 = \operatorname{tr}(\Sigma^2)$, and that $\mathsf{E}(X_1'X_1) = \operatorname{tr}(\Sigma)$.
Thus, we obtain that
\begin{eqnarray*}
&& \mathsf{E}\bigl[\bigl(X_1'X_2
\bigr)^2 - \bigl(X_1'X_1\bigr)
- \bigl(X_2'X_2\bigr)\bigr] \bigl[
\bigl(X_1'X_3\bigr)^2 -
\bigl(X_1'X_1\bigr) -
\bigl(X_3'X_3\bigr)\bigr]
\\
&&\quad  = \operatorname{tr}^2\bigl(\Sigma^2\bigr) + 2
\operatorname{tr}\bigl(\Sigma^4\bigr) - 4\operatorname{tr}\bigl(
\Sigma^3\bigr) - 4\operatorname{tr}\bigl(\Sigma^2\bigr)
\operatorname{tr}(\Sigma )+2\operatorname{tr}\bigl(\Sigma^2\bigr)+4
\operatorname{tr}^2(\Sigma).
\end{eqnarray*}
Noting that $\mathsf{E}[(X_1'X_2)^2 - (X_1'X_1) - (X_2'X_2)] =
\operatorname{tr}(\Sigma^2)-2\operatorname{tr}(\Sigma)$, we obtain the claim.
\end{pf}

\begin{pf*}{Proof of \eqref{eq:var-T_n}}
With Lemma \ref{lemma:var-cov-h}, we have
\begin{eqnarray*}
&&\operatorname{\mathsf{Var}} \biggl(\sum_{1\leq i<j\leq n}
h(X_i,X_j) \biggr)\\
 &&\quad  = \sum
_{1\leq
i<j\leq n} \operatorname{\mathsf{Var}} \bigl(h(X_i,X_j)
\bigr) + 2
\mathop{\sum_{1\leq i<j, i'<j'\leq n}}_{i=i' \mathrm{~or~} j=j'}
 \operatorname{\mathsf{Cov}}
\bigl( h(X_i,X_j), h(X_{i'},
X_{j'}) \bigr)
\\
&&\quad  = \frac{n(n-1)}{2}\operatorname{\mathsf{Var}} \bigl(h(X_1,X_2)
\bigr) + 2\frac
{n(n-1)}{2}(n-2)\operatorname{\mathsf{Cov}} \bigl(
h(X_1,X_2), h(X_1, X_3) \bigr)
\\
&&\quad  = n(n-1) \bigl[\operatorname{tr}^2\bigl(\Sigma^2\bigr)
+ \operatorname {tr}\bigl(\Sigma^4\bigr) \bigr] +
2n(n-1)^2\operatorname{tr} \bigl(\Sigma^2(
\Sigma-I)^2 \bigr).
\end{eqnarray*}
Multiplying both sides with $4n^{-2}(n-1)^{-2}$, we obtain \eqref{eq:var-T_n}.
\end{pf*}

%s6.3 ###
\subsection{\texorpdfstring{Proof details for Proposition \protect\ref{prop:T-n-MGCLT}}{Proof details for Proposition 3}}
\label{app:prop-rate-detail}

%s6.3.1 ###
\subsubsection{\texorpdfstring{Proof of \protect\eqref{eq:T_n-MG-diff}}{Proof of (12)}}

First of all, we give a formal proof of the representation \eqref
{eq:T_n-MG-diff}.

\begin{pf*}{Proof of \eqref{eq:T_n-MG-diff}}
The computations made in \cite{chzh10}, Appendix, are handy for the
proof here. Indeed, we have
\begin{eqnarray}
\label{eq:u-v} %
 u_{nk} & =& (
\mathsf{E}_k - \mathsf{E}_{k-1}) \biggl[
\frac{1}{n}X_k'X_k \biggr] =
\frac{1}{n} \bigl[X_k'X_k -
\operatorname{tr}(\Sigma) \bigr],
\\
v_{nk} & =& (\mathsf{E}_k - \mathsf{E}_{k-1})
\biggl[ \frac
{2}{n(n-1)}\sum_{i\neq
k}
\bigl(X_i'X_k\bigr)^2 \biggr]\nonumber
\\
& =& \frac{2}{n(n-1)}\bigl[ X_k'Q_{k-1}X_k
- \operatorname {tr}(Q_{k-1}\Sigma) \bigr] + \frac
{2}{n}\bigl[
X_k'\Sigma X_k - \operatorname{tr}\bigl(
\Sigma^2\bigr) \bigr],\nonumber
\end{eqnarray}
where $Q_{k-1} = \sum_{i=1}^{k-1}(X_iX_i' - \Sigma)$. Noting that
$D_{nk} = v_{nk}-2u_{nk}$, we obtain \eqref{eq:T_n-MG-diff}.
\end{pf*}

%s6.3.2 ###
\subsubsection{\texorpdfstring{Proof of \protect\eqref{eq:sigma_nk}}{Proof of (19)}}

To calculate $\sigma_{nk}^2$, we note that
\[
\sigma_{nk}^2 = \mathsf{E}_{k-1}
\bigl[D_{nk}^2\bigr] = 4\mathsf {E}_{k-1}
\bigl[u_{nk}^2\bigr] - 4\mathsf{E}_{k-1}[u_{nk}v_{nk}]
+ \mathsf{E}_{k-1}\bigl[v_{nk}^2\bigr].
\]
Thus, \eqref{eq:sigma_nk} is immediate with the following lemma.

\begin{lemma}
\label{lemma:eq-u-v-squared}
For $u_{nk}, v_{nk}$ defined as in \eqref{eq:u-v}, we have
% \begin{equation}
% \label{eq:u-v-squared}
%
\begin{eqnarray*}
\mathsf{E}_{k-1}\bigl[u_{nk}^2\bigr] & =&
\frac{2}{n^2}\operatorname{tr}\bigl(\Sigma^2\bigr),
\\
\mathsf{E}_{k-1}[u_{nk}v_{nk}] & =&
\frac{4}{n^2(n-1)}\operatorname {tr}\bigl(Q_{k-1}\Sigma^2
\bigr) + \frac{4}{n^2}\operatorname{tr}\bigl(\Sigma^3\bigr),
\\
\mathsf{E}_{k-1}\bigl[v_{nk}^2\bigr] & =&
\frac{8}{n^2(n-1)^2}\operatorname {tr}(Q_{k-1}\Sigma Q_{k-1}
\Sigma)\\
&&{} + \frac{16}{n^2(n-1)}\operatorname{tr}\bigl(Q_{k-1}
\Sigma^3\bigr) + \frac
{8}{n^2}\operatorname{tr}\bigl(
\Sigma^4\bigr).
\end{eqnarray*}
%
% \end{equation}
\end{lemma}

\begin{pf}
First, we have
\[
\mathsf{E}_{k-1}\bigl[u_{nk}^2\bigr] =
\frac{1}{n^2}\operatorname{\mathsf {Var}}\bigl(X_k'X_k
\bigr) = \frac{2}{n^2}\operatorname{tr} \bigl(\Sigma^2\bigr).
\]
Next, we have from \eqref{eq:cross-2nd-moment} that
\begin{eqnarray*}
\mathsf{E}_{k-1}[u_{nk}v_{nk}] & =&
\frac{2}{n^2(n-1)}\mathsf {E}_{k-1} \bigl[X_k'X_k
- \operatorname{tr}(\Sigma) \bigr]X_k'Q_{k-1}X_k
+ \frac{2}{n^2} \bigl[X_k'X_k -
\operatorname{tr}(\Sigma) \bigr]X_k'\Sigma
X_k
\\
& = &\frac{2}{n^2(n-1)} \bigl[\operatorname{tr}(\Sigma)\operatorname
{tr}(Q_{k-1}\Sigma) + 2\operatorname{tr} \bigl(Q_{k-1}
\Sigma^2\bigr) - \operatorname{tr}(\Sigma)\operatorname
{tr}(Q_{k-1}\Sigma) \bigr]
\\
&&{}  + \frac{2}{n^2} \bigl[\operatorname{tr}(\Sigma)\operatorname {tr}
\bigl(\Sigma^2\bigr) + 2\operatorname{tr} \bigl(\Sigma^3
\bigr) - \operatorname{tr}(\Sigma)\operatorname{tr}\bigl(\Sigma^2
\bigr) \bigr]
\\
& = &\frac{4}{n^2(n-1)}\operatorname{tr}\bigl(Q_{k-1}
\Sigma^2\bigr)+\frac
{4}{n^2}\operatorname{tr} \bigl(
\Sigma^3\bigr).
\end{eqnarray*}
Finally, we have
\begin{eqnarray*}
\mathsf{E}_{k-1}\bigl[v_{nk}^2\bigr] & =&
\frac{4}{n^2(n-1)^2}\mathsf {E}_{k-1} \bigl[X_k'Q_{k-1}X_k
- \operatorname{tr}(Q_{k-1}\Sigma) \bigr]X_k'
Q_{k-1}X_k
\\
& &{} + \frac{8}{n^2(n-1)}\mathsf{E}_{k-1} \bigl[X_k'
\Sigma X_k - \operatorname{tr} \bigl(\Sigma^2\bigr)
\bigr]X_k' Q_{k-1}X_k
\\
&&{}  + \frac{4}{n^2}\mathsf{E}_{k-1} \bigl[X_k'
\Sigma X_k - \operatorname{tr}\bigl(\Sigma^2\bigr)
\bigr]X_k' \Sigma X_k.
\end{eqnarray*}
Note that
\begin{eqnarray*}
\mathsf{E}_{k-1} \bigl[X_k'Q_{k-1}X_k
X_k'Q_{k-1}X_k \bigr] & =&
\mathsf {E}_{k-1} \bigl[Z_k'
\Gamma'Q_{k-1}\Gamma Z_k Z_k'
\Gamma'Q_{k-1}\Gamma Z_k \bigr]
\\
& = &\operatorname{tr}^2\bigl(\Gamma'Q_{k-1}
\Gamma\bigr) + 2\operatorname {tr}\bigl(\Gamma'Q_{k-1}
\Gamma\Gamma' Q_{k-1}\Gamma\bigr)
\\
& =& 2\operatorname{tr}(Q_{k-1}\Sigma Q_{k-1}\Sigma) +
\operatorname {tr}^2(Q_{k-1}\Sigma),
\\
\mathsf{E} \bigl[X_k'Q_{k-1}X_k
X_k'\Sigma X_k \bigr] & =& \mathsf
{E}_{k-1} \bigl[Z_k'\Gamma'Q_{k-1}
\Gamma Z_k Z_k'\Gamma'\Sigma
\Gamma Z_k \bigr]
\\
& =& \operatorname{tr}\bigl(\Gamma'Q_{k-1}\Gamma\bigr)
\operatorname{tr}\bigl(\Gamma'\Sigma\Gamma\bigr) + 2
\operatorname{tr} \bigl(\Gamma'Q_{k-1}\Gamma
\Gamma'\Sigma\Gamma\bigr)
\\
& =& 2\operatorname{tr}\bigl(Q_{k-1}\Sigma^3\bigr) +
\operatorname {tr}(Q_{k-1}\Sigma)\operatorname{tr}\bigl(
\Sigma^2\bigr),
\\
\mathsf{E}_{k-1} \bigl[X_k'\Sigma
X_k X_k'\Sigma X_k \bigr] &
= &\mathsf{E}_{k-1} \bigl[Z_k'
\Gamma'\Sigma\Gamma Z_k Z_k'
\Gamma'\Sigma\Gamma Z_k \bigr]
\\
& = &2\operatorname{tr}\bigl(\Gamma'\Sigma\Gamma
\Gamma'\Sigma\Gamma\bigr) + \operatorname{tr}^2\bigl(
\Gamma'\Sigma\Gamma\bigr)
\\
& =& 2\operatorname{tr}\bigl(\Sigma^4\bigr) + \operatorname{tr}^2
\bigl(\Sigma^2\bigr).
\end{eqnarray*}
Collecting terms, we complete the proof.
\end{pf}

%s6.3.3 ###
\subsubsection{\texorpdfstring{Proof of Lemma \protect\ref{lemma:tr-M-k-1-square}}{Proof of Lemma 3}}

Finally, we shall complete the proof of Lemma \ref{lemma:tr-M-k-1-square}.

Recall that $M_{k-1} = A\sum_{i=1}^{k-1}(Z_iZ_i'-I)A$, and so
\begin{eqnarray}
\label{eq:tr-M-k-1-sq} \operatorname{tr}\bigl(M_{k-1}^2
\bigr) & =& \sum_{i=1}^{k-1}\operatorname {tr}
\bigl(A\bigl(Z_iZ_i'-I
\bigr)A^2\bigl(Z_iZ_i'-I
\bigr)A\bigr)
\nonumber\\[-8pt]\\[-8pt]
&&{}+ 2\sum_{i=1}^{k-2}\sum
_{j=i+1}^{k-1}\operatorname {tr}\bigl(A
\bigl(Z_iZ_i'-I\bigr)A^2
\bigl(Z_jZ_j'-I\bigr)A\bigr).\nonumber
\end{eqnarray}
For any fixed $i$, we have
\begin{eqnarray*}
\mathsf{E}\bigl[\operatorname{tr}\bigl(A\bigl(Z_iZ_i'-I
\bigr)A^2\bigl(Z_iZ_i'-I
\bigr)A\bigr)\bigr] & =& \mathsf{E}\bigl[\operatorname{tr}\bigl(AZ_iZ_i'A^2
Z_iZ_i'A\bigr)\bigr] - 2\mathsf{E}\bigl[
\operatorname{tr}\bigl(A^2 Z_i Z_i'
A^2\bigr)\bigr] + \operatorname{tr}\bigl(A^4\bigr)
\\
& =& \mathsf{E}\bigl[\bigl(Z_i'A^2
Z_i\bigr)^2\bigr] - 2\mathsf{E}\bigl[Z_i'A^4
Z_i\bigr] + \operatorname{tr}\bigl(A^4\bigr)
\\
& =& 2\operatorname{tr}\bigl(A^4\bigr) + \operatorname{tr}^2
\bigl(A^2\bigr) - 2\operatorname{tr}\bigl(A^4\bigr) +
\operatorname{tr}\bigl(A^4\bigr)
\\
& =& \operatorname{tr}^2\bigl(\Sigma^2\bigr) +
\operatorname{tr}\bigl(\Sigma^4\bigr).
\end{eqnarray*}
On the other hand, for any $i\neq j$, we have
\[
\mathsf{E}\bigl[\operatorname{tr}\bigl(A\bigl(Z_iZ_i'-I
\bigr)A^2\bigl(Z_jZ_j'-I
\bigr)A\bigr)\bigr]  = 0.
\]
In addition, we note that the terms in \eqref{eq:tr-M-k-1-sq} are all
uncorrelated. Therefore, we obtain that
\[
\mathsf{E}\bigl[\operatorname{tr}\bigl(M_{k-1}^2\bigr)
\bigr]  = (k-1)\bigl[\operatorname {tr}^2\bigl(\Sigma^2
\bigr) + \operatorname{tr}\bigl(\Sigma^4\bigr)\bigr].
\]
Moreover, we have
\begin{eqnarray}
\label{eq:V1+V2} %
 \operatorname{\mathsf{Var}}\bigl[
\operatorname{tr}\bigl(M_{k-1}^2\bigr)\bigr] & =& \sum
_{i=1}^{k-1}\operatorname{\mathsf{Var}}\bigl[
\operatorname{tr} \bigl(A\bigl(Z_iZ_i'-I
\bigr)A^2\bigl(Z_iZ_i'-I
\bigr)A\bigr)\bigr]\nonumber
\\
&&{}  + 2\sum_{i=1}^{k-2}\sum
_{j=i+1}^{k-1}\operatorname{\mathsf {Var}}\bigl[
\operatorname{tr} \bigl(A\bigl(Z_iZ_i'-I
\bigr)A^2\bigl(Z_jZ_j'-I
\bigr)A\bigr)\bigr]
\\
& =& (k-1)V_1 + 2(k-1) (k-2)V_2. \nonumber
\end{eqnarray}
Here, $V_1 = \operatorname{\mathsf{Var}}[\operatorname
{tr}(A(Z_1Z_1'-I)A^2(Z_1Z_1'-I)A)]$ and $V_2 = \operatorname{\mathsf{Var}}
[\operatorname{tr}(A(Z_1Z_1'-I)A^2(Z_2Z_2'-I)A)]$.

Consider $V_1$ first, for which we have the decomposition
\begin{eqnarray*}
V_1 & =& \operatorname{\mathsf{Var}}\bigl[\bigl(Z_1'A^2
Z_1\bigr)^2 - 2Z_1'A^4
Z_1\bigr]
\\
& =& \operatorname{\mathsf{Var}}\bigl[\bigl(Z_1'A^2
Z_1\bigr)^2\bigr] + 4\operatorname {\mathsf{Var}}
\bigl(Z_1'A^4 Z_1\bigr) - 4
\operatorname{\mathsf{Cov}}\bigl[\bigl(Z_1'A^2
Z_1\bigr)^2, Z_1'A^4
Z_1\bigr].
\end{eqnarray*}

To calculate $\operatorname{\mathsf{Var}}[(Z_1'A^2 Z_1)^2]$, we note
that the eigenvalues of
$A^2$ are $\lambda_1^2\geq\cdots\geq\lambda_p^2$, where $\lambda_1\geq\cdots\geq\lambda_p$ are the eigenvalues of $\Sigma$. Let
$U_j\stackrel{\mathrm{i.i.d.}}{\sim}N(0,1)$, for $j=1,\dots, p$. By the moment generating
function of $\chi^2_{(1)}$ distribution, we have $\mathsf{E}[U_j^2] = 1$,
$\mathsf{E}[U_j^4] = 3$, $\mathsf{E}[U_j^6] = 15$, and $\mathsf
{E}[U_j^8] = 105$. Then,
we obtain that
\begin{eqnarray*}
\mathsf{E}\bigl(Z_i'A^2 Z_i
\bigr)^4 & =& \mathsf{E} \Biggl[ \Biggl(\sum
_{j=1}^p\lambda_j^2
U_j^2 \Biggr)^4 \Biggr]
\\
& =& \mathsf{E} \Biggl[ \sum_{j=1}^p
\lambda_j^8 U_j^8 + 4\sum
_{j\neq l}^p \lambda_j^6
\lambda_l^2 U_j^6
U_l^2 + \frac{6}{2!}\sum
_{j\neq l}^p \lambda_j^4
\lambda_l^4 U_j^4
U_l^4
\\
&&\hspace*{12pt}{}   + \frac{12}{2!}\sum_{j\neq l\neq m}^p
\lambda_j^4\lambda_l^2
\lambda_m^2 U_j^4
U_l^2 U_m^2 +
\frac{24}{4!}\sum_{j\neq l\neq m\neq r}^p
\lambda_j^2\lambda_l^2
\lambda_m^2\lambda_r^2
U_j^2 U_l^2
U_m^2 U_r^2 \Biggr]
\\
& =& 105\sum_{j=1}^p \lambda_j^8
+ 15\cdot4\sum_{j\neq l}^p
\lambda_j^6\lambda_l^2 + 9\cdot
\frac{6}{2!}\sum_{j\neq l}^p
\lambda_j^4\lambda_l^4
\\
&&{}  + 3\cdot\frac{12}{2!}\sum_{j\neq l\neq m}^p
\lambda_j^4\lambda_l^2
\lambda_m^2 + \frac{24}{4!}\sum
_{j\neq l\neq m\neq
r}^p \lambda_j^2
\lambda_l^2\lambda_m^2
\lambda_r^2
\\
& = &\operatorname{tr}^4\bigl(\Sigma^2\bigr) + 12
\operatorname{tr}\bigl(\Sigma^4\bigr)\operatorname{tr}^2
\bigl(\Sigma^2\bigr) + 12\operatorname{tr}^2\bigl(
\Sigma^4\bigr) + 32\operatorname{tr}\bigl(\Sigma^6\bigr)
\operatorname {tr}\bigl(\Sigma^2\bigr) + 48\operatorname{tr}\bigl(
\Sigma^8\bigr).
\end{eqnarray*}
Observing that $\mathsf{E}[(Z_1'A^2 Z_1)^2] = 2\operatorname
{tr}(\Sigma^4)+\operatorname{tr}^2(\Sigma^2)$, we get\vspace*{1pt}
\[
\operatorname{\mathsf{Var}}\bigl[\bigl(Z_1'A^2
Z_1\bigr)^2\bigr] = 8\operatorname {tr}\bigl(
\Sigma^4\bigr)\operatorname{tr}^2\bigl(
\Sigma^2\bigr) + 8\operatorname{tr}^2\bigl(
\Sigma^4\bigr) + 32\operatorname{tr}\bigl(\Sigma^6\bigr)
\operatorname {tr}\bigl(\Sigma^2\bigr) + 48\operatorname{tr}\bigl(
\Sigma^8\bigr).
\]

Next, we compute $\operatorname{\mathsf{Var}}(Z_1'A^4 Z_1)$, for
which we have\vspace*{1pt}
\begin{eqnarray*}
\mathsf{E}\bigl[\bigl(Z_1'A^4
Z_1\bigr)^2\bigr] &=& 2\operatorname{tr}\bigl(
\Sigma^8\bigr) + \operatorname{tr}^2\bigl(
\Sigma^4\bigr),\\
  \mathsf{E}\bigl[Z_1'A^4
Z_1\bigr] &=& \operatorname{tr}\bigl(\Sigma^4\bigr).
\end{eqnarray*}
Therefore, we get $\operatorname{\mathsf{Var}}(Z_1'A^4 Z_1) =
2\operatorname{tr}(\Sigma^8)$.

Now switch to $\operatorname{\mathsf{Cov}}[(Z_1'A^2 Z_1)^2, Z_1'A^4
Z_1]$. We note that\vspace*{1pt}
\begin{eqnarray*}
&& \mathsf{E}\bigl[\bigl(Z_1'A^2
Z_1\bigr)^2 Z_1'
A^4 Z_1\bigr]
\\
&&\quad =  \mathsf{E} \Biggl[ \Biggl( \sum_{j=1}^p
\lambda_j^4 U_j^4 + \sum
_{j\neq
l}^p \lambda_j^2
\lambda_l^2 U_j^2
U_l^2 \Biggr) \sum_{j=1}^p
\lambda_j^4 U_j^2 \Biggr]
\\
&&\quad = \mathsf{E} \Biggl[ \sum_{j=1}^p
\lambda_j^8 U_j^6 + \sum
_{j\neq l}^p \lambda_j^4
\lambda_l^4 U_j^4
U_l^2 + 2\sum_{j\neq l}^p
\lambda_j^6\lambda_l^2
U_j^4 U_l^2 + \sum
_{j\neq l\neq m}^p \lambda_j^4
\lambda_l^2\lambda_m^2
U_j^2 U_l^2
U_m^2 \Biggr]
\\
&&\quad = 15\sum_{j=1}^p\lambda_j^8
+ 3\sum_{j\neq l}^p \lambda_j^4
\lambda_l^4 + 6\sum_{j\neq l}^p
\lambda_j^6 \lambda_l^2 + \sum
_{j\neq l\neq m}^p \lambda_j^4
\lambda_l^2\lambda_m^2
\\
&&\quad = \operatorname{tr}\bigl(\Sigma^4\bigr)\operatorname{tr}^2
\bigl(\Sigma^2\bigr) + 4 \operatorname{tr}\bigl(\Sigma^6
\bigr)\operatorname{tr}\bigl(\Sigma^2\bigr) + 2\operatorname{tr}^2
\bigl(\Sigma^4\bigr) + 8 \operatorname{tr}\bigl(\Sigma^8
\bigr).
\end{eqnarray*}
By previous expression for $\mathsf{E}[(Z_1'A^2 Z_1)^2]$ and $\mathsf
{E}[Z_1'A^4
Z_1]$, we obtain that\vspace*{1pt}
\[
\operatorname{\mathsf{Cov}}\bigl[\bigl(Z_1'A^2
Z_1\bigr)^2, Z_1'A^4
Z_1\bigr] = 4\operatorname{tr}\bigl(\Sigma^6\bigr)
\operatorname{tr}\bigl(\Sigma^2\bigr) + 8 \operatorname{tr}\bigl(
\Sigma^8\bigr).
\]\eject

Finally, we obtain that
\begin{equation}
\label{eq:V1} V_1 = 8\operatorname{tr}\bigl(\Sigma^4
\bigr)\operatorname{tr}^2\bigl(\Sigma^2\bigr) + 8
\operatorname{tr}^2\bigl(\Sigma^4\bigr) + 16
\operatorname{tr} \bigl(\Sigma^6\bigr)\operatorname{tr}\bigl(
\Sigma^2\bigr) + 24\operatorname{tr}\bigl(\Sigma^8\bigr).
\end{equation}

Switch to the calculation of $V_2$. We first note that
\begin{eqnarray*}
V_2 & =& \operatorname{\mathsf{Var}}\bigl[\operatorname {tr}\bigl(A
\bigl(Z_1Z_1'-I\bigr)A^2
\bigl(Z_2Z_2'-I\bigr)A\bigr)\bigr]
\\
& =& \operatorname{\mathsf{Var}}\bigl[\operatorname{tr}\bigl(A Z_1
Z_1' A^2 Z_1
Z_1' A\bigr) - \operatorname{tr}\bigl(A^2
Z_1 Z_1' A^2\bigr) -
\operatorname{tr}\bigl(A^2 Z_2 Z_2'
A^2\bigr)\bigr]
\\
& =& \operatorname{\mathsf{Var}}\bigl[\bigl(Z_1'A^2
Z_2\bigr)^2\bigr] + 2\operatorname {\mathsf{Var}}
\bigl(Z_1'A^4 Z_1\bigr) - 2
\operatorname{\mathsf{Cov}}\bigl[\bigl(Z_1'A^2
Z_2\bigr)^2, Z_1'A^4
Z_1\bigr].
\end{eqnarray*}
Note that $\mathsf{E}[(Z_1' A^2 Z_2)^4] = 3\operatorname{tr}^2(\Sigma^4) + 6\operatorname{tr}(\Sigma^8)$, and $\mathsf{E}[(Z_1'A^2 Z_2)^2] = \operatorname{tr}(\Sigma^4)$. We then get
\[
\operatorname{\mathsf{Var}}\bigl[\bigl(Z_1'
A^2 Z_2\bigr)^2\bigr] = 2\operatorname
{tr}^2\bigl(\Sigma^4\bigr) + 6\operatorname{tr}\bigl(
\Sigma^8\bigr).
\]
In addition, previous calculation gives $\operatorname{\mathsf
{Var}}(Z_1'A^4 Z_1) = 2\operatorname{tr}
(\Sigma^8)$.
Then for $\operatorname{\mathsf{Cov}}[(Z_1'A^2 Z_2)^2,\allowbreak  Z_1'A^4 Z_1]$,
we have
\begin{eqnarray*}
\mathsf{E}\bigl[\bigl(Z_1' A^2
Z_1\bigr)^2 Z_1'
A^4 Z_1\bigr] & =& \mathsf{E} \Biggl[ \Biggl( \sum
_{j=1}^p \lambda_j^2
U_j V_j \Biggr)^2 \sum
_{j=1}^p \lambda_j^4
U_j^2 \Biggr]
\\
& =& \mathsf{E} \Biggl[ \Biggl(\sum_{j=1}^p
\lambda_j^4 U_j^2
V_j^2 \Biggr) \Biggl(\sum_{j=1}^p
\lambda_j^4 U_j^2 \Biggr)
\Biggr]
\\
& =& \mathsf{E} \Biggl[ \sum_{j=1}^p
\lambda_j^8 U_j^4
V_j^2 + \sum_{j\neq
l}^p
\lambda_j^4 \lambda_l^4
U_j^2 V_j^2
U_l^2 \Biggr]
\\
& =& 3\sum_{j=1}^p \lambda_j^8
+ \sum_{j\neq l}^p \lambda_j^4
\lambda_l^4
\\
& =& \operatorname{tr}^2\bigl(\Sigma^4\bigr) + 2
\operatorname{tr}\bigl(\Sigma^8\bigr).
\end{eqnarray*}
This leads to the conclusion that $\operatorname{\mathsf
{Cov}}[(Z_1'A^2 Z_2)^2, Z_1'A^4 Z_1]
= 2\operatorname{tr}(\Sigma^8)$, and so
\begin{equation}
\label{eq:V2} V_2 = 2\operatorname{tr}^2\bigl(
\Sigma^4\bigr) + 6\operatorname{tr}\bigl(\Sigma^8\bigr).
\end{equation}
Replacing $V_1$ and $V_2$ in \eqref{eq:V1+V2} by \eqref{eq:V1} and
\eqref{eq:V2}, we obtain the claimed formula for $\operatorname
{\mathsf{Var}}[\operatorname{tr}(M_{k-1}^2)]$.
\end{appendix}
\iffalse
\fi
\section*{Acknowledgements}
The research of Tony Cai was supported in part by NSF FRG Grant DMS-08-54973.
The research of Zongming Ma was supported in part by the Dean's
Research Fund of The Wharton School.

%suskaldyti doi

% imsref loaded by imikolaityte, 2012-08-20 09:40:09
%

\printhistory


\begin{thebibliography}{19}
% BibTex style file: bej.bst, 2011-10-13
% Default style options (sort=1,type=number).
% Used options (sort=1,type=number).

%b1 ###
\bibitem{ande03}
%
\begin{bbook}[mr]
\bauthor{\bsnm{Anderson},~\bfnm{T.~W.}\binits{T.W.}}
(\byear{2003}).
\btitle{An Introduction to Multivariate Statistical Analysis},
\bedition{3rd} ed.
\bseries{Wiley Series in Probability and Statistics}.
\baddress{Hoboken, NJ}: \bpublisher{Wiley-Interscience [John Wiley \& Sons]}.
\bid{mr={1990662}}
\bptok{imsref}%
\end{bbook}
%
\endbibitem

%b2 ###
\bibitem{baji09}
%
\begin{barticle}[mr]
\bauthor{\bsnm{Bai},~\bfnm{Zhidong}\binits{Z.}},
\bauthor{\bsnm{Jiang},~\bfnm{Dandan}\binits{D.}},
\bauthor{\bsnm{Yao},~\bfnm{Jian-Feng}\binits{J.F.}} \AND
\bauthor{\bsnm{Zheng},~\bfnm{Shurong}\binits{S.}}
(\byear{2009}).
\btitle{Corrections to {LRT} on large-dimensional covariance matrix by {RMT}}.
\bjournal{Ann. Statist.}
\bvolume{37}
\bpages{3822--3840}.
\bid{doi={10.1214/09-AOS694}, issn={0090-5364}, mr={2572444}}
\bptok{imsref}%
\end{barticle}
%
\endbibitem

%b3 ###
\bibitem{bide05}
%
\begin{barticle}[mr]
\bauthor{\bsnm{Birke},~\bfnm{Melanie}\binits{M.}} \AND
\bauthor{\bsnm{Dette},~\bfnm{Holger}\binits{H.}}
(\byear{2005}).
\btitle{A note on testing the covariance matrix for large dimension}.
\bjournal{Statist. Probab. Lett.}
\bvolume{74}
\bpages{281--289}.
\bid{doi={10.1016/j.spl.2005.04.051}, issn={0167-7152}, mr={2189467}}
\bptok{imsref}%
\end{barticle}
%
\endbibitem

%b5 ###
\bibitem{caji11}
%
\begin{barticle}[mr]
\bauthor{\bsnm{Cai},~\bfnm{T.~Tony}\binits{T.T.}} \AND
\bauthor{\bsnm{Jiang},~\bfnm{Tiefeng}\binits{T.}}
(\byear{2011}).
\btitle{Limiting laws of coherence of random matrices with
applications to
testing covariance structure and construction of compressed sensing
matrices}.
\bjournal{Ann. Statist.}
\bvolume{39}
\bpages{1496--1525}.
\bid{doi={10.1214/11-AOS879}, issn={0090-5364}, mr={2850210}}
\bptok{imsref}%
\end{barticle}
%
\endbibitem

%b4 ###
\bibitem{cai11}
%
\begin{bmisc}[auto:STB|2012/08/14|15:18:37]
\bauthor{\bsnm{Cai},~\bfnm{T.~T.}\binits{T.T.}},
\bauthor{\bsnm{Liu},~\bfnm{W.}\binits{W.}} \AND
\bauthor{\bsnm{Xia},~\bfnm{Y.}\binits{Y.}}
(\byear{2011}).
\bhowpublished{Two-sample covariance matrix
testing and support recovery. Technical report}.
\bptok{imsref}%
\end{bmisc}
%
\endbibitem



%b6 ###
\bibitem{chenli12}
%
\begin{barticle}[auto:STB|2012/08/14|15:18:37]
\bauthor{\bsnm{Chen},~\bfnm{S.~X.}\binits{S.X.}}\AND
\bauthor{\bsnm{Li},~\bfnm{J.}\binits{J.}}
(\byear{2012}).
\btitle{Two sample tests for high dimensional
covariance matrices}.
\bjournal{Ann. Statist.}
\bvolume{40}
\bpages{908--940}.
\bid{mr={2985938}}
\bptok{imsref}%
\end{barticle}
%
\endbibitem

%b7 ###
\bibitem{chzh10}
%
\begin{barticle}[mr]
\bauthor{\bsnm{Chen},~\bfnm{Song~Xi}\binits{S.X.}},
\bauthor{\bsnm{Zhang},~\bfnm{Li-Xin}\binits{L.X.}} \AND
\bauthor{\bsnm{Zhong},~\bfnm{Ping-Shou}\binits{P.S.}}
(\byear{2010}).
\btitle{Tests for high-dimensional covariance matrices}.
\bjournal{J.~Amer. Statist. Assoc.}
\bvolume{105}
\bpages{810--819}.
\bid{doi={10.1198/jasa.2010.tm09560}, issn={0162-1459}, mr={2724863}}
\bptok{imsref}%
\end{barticle}
%
\endbibitem

%b8 ###
\bibitem{hebr70}
%
\begin{barticle}[mr]
\bauthor{\bsnm{Heyde},~\bfnm{C.~C.}\binits{C.C.}} \AND
\bauthor{\bsnm{Brown},~\bfnm{B.~M.}\binits{B.M.}}
(\byear{1970}).
\btitle{On the departure from normality of a certain class of martingales}.
\bjournal{Ann. Math. Statist.}
\bvolume{41}
\bpages{2161--2165}.
\bid{issn={0003-4851}, mr={0293702}}
\bptok{imsref}%
\end{barticle}
%
\endbibitem

%b9 ###
\bibitem{hoef63}
%
\begin{barticle}[mr]
\bauthor{\bsnm{Hoeffding},~\bfnm{Wassily}\binits{W.}}
(\byear{1963}).
\btitle{Probability inequalities for sums of bounded random variables}.
\bjournal{J.~Amer. Statist. Assoc.}
\bvolume{58}
\bpages{13--30}.
\bid{issn={0162-1459}, mr={0144363}}
\bptok{imsref}%
\end{barticle}
%
\endbibitem


%b10 ###
\bibitem{jiji10}
%
\begin{barticle}[auto:STB|2012/08/14|15:18:37]
\bauthor{\bsnm{Jiang},~\bfnm{D.}\binits{D.}},
\bauthor{\bsnm{Jiang},~\bfnm{T.}\binits{T.}} \AND
\bauthor{\bsnm{Yang},~\bfnm{F.}\binits{F.}}
(\byear{2012}).
\btitle{Likelihood ratio tests for
covariance matrices of high-dimensional normal distributions}.
\bjournal{J. Statist. Plann. Inference}
\bvolume{142}
\bpages{2241--2256}.
\bid{mr={291184}}
\bptok{imsref}%
\end{barticle}
%
\endbibitem

%b11 ###
\bibitem{john01}
%
\begin{barticle}[mr]
\bauthor{\bsnm{Johnstone},~\bfnm{Iain~M.}\binits{I.M.}}
(\byear{2001}).
\btitle{On the distribution of the largest eigenvalue in principal components
analysis}.
\bjournal{Ann. Statist.}
\bvolume{29}
\bpages{295--327}.
\bid{doi={10.1214/aos/1009210544}, issn={0090-5364}, mr={1863961}}
\bptok{imsref}%
\end{barticle}
%
\endbibitem

%b12 ###
\bibitem{lewo02}
%
\begin{barticle}[mr]
\bauthor{\bsnm{Ledoit},~\bfnm{Olivier}\binits{O.}} \AND
\bauthor{\bsnm{Wolf},~\bfnm{Michael}\binits{M.}}
(\byear{2002}).
\btitle{Some hypothesis tests for the covariance matrix when the
dimension is
large compared to the sample size}.
\bjournal{Ann. Statist.}
\bvolume{30}
\bpages{1081--1102}.
\bid{doi={10.1214/aos/1031689018}, issn={0090-5364}, mr={1926169}}
\bptok{imsref}%
\end{barticle}
%
\endbibitem

%b13 ###
\bibitem{muir82}
%
\begin{bbook}[mr]
\bauthor{\bsnm{Muirhead},~\bfnm{Robb~J.}\binits{R.J.}}
(\byear{1982}).
\btitle{Aspects of Multivariate Statistical Theory}.
\bseries{Wiley Series in Probability and Mathematical Statistics}.
\baddress{New York}: \bpublisher{Wiley}.
\bid{mr={0652932}}
\bptok{imsref}%
\end{bbook}
%
\endbibitem

%b14 ###
\bibitem{naga73}
%
\begin{barticle}[mr]
\bauthor{\bsnm{Nagao},~\bfnm{Hisao}\binits{H.}}
(\byear{1973}).
\btitle{On some test criteria for covariance matrix}.
\bjournal{Ann. Statist.}
\bvolume{1}
\bpages{700--709}.
\bid{issn={0090-5364}, mr={0339405}}
\bptok{imsref}%
\end{barticle}
%
\endbibitem

%b15 ###
\bibitem{onmo11}
%
\begin{bmisc}[auto:STB|2012/08/14|15:18:37]
\bauthor{\bsnm{Onatski},~\bfnm{A.}\binits{A.}},
\bauthor{\bsnm{Moreira},~\bfnm{M.J.}\binits{M.J.}} \AND
\bauthor{\bsnm{Hallin},~\bfnm{M.}\binits{M.}}
(\byear{2011}).
\bhowpublished{Asymptotic power of
sphericity tests for high-dimensional data. Available at
\texttt{\href{http://www.econ.cam.ac.uk/people/faculty/ao319/pubs/WPOnatskiMoreira.pdf}{http://www.econ.cam.ac.uk/people/}
\href{http://www.econ.cam.ac.uk/people/faculty/ao319/pubs/WPOnatskiMoreira.pdf}{faculty/ao319/pubs/WPOnatskiMoreira.pdf}}}.
\bptok{imsref}%
\end{bmisc}
%
\endbibitem

%b16 ###
\bibitem{roy57}
%
\begin{bbook}[mr]
\bauthor{\bsnm{Roy},~\bfnm{S.~N.}\binits{S.N.}}
(\byear{1957}).
\btitle{Some Aspects of Multivariate Analysis}.
\baddress{New York}: \bpublisher{Wiley}.
\bid{mr={0092296}}
\bptok{imsref}%
\end{bbook}
%
\endbibitem

%b17 ###
\bibitem{sriv05}
%
\begin{barticle}[mr]
\bauthor{\bsnm{Srivastava},~\bfnm{Muni~S.}\binits{M.S.}}
(\byear{2005}).
\btitle{Some tests concerning the covariance matrix in high
dimensional data}.
\bjournal{J.~Japan Statist. Soc.}
\bvolume{35}
\bpages{251--272}.
\bid{issn={0389-5602}, mr={2328427}}
\bptok{imsref}%
\end{barticle}
%
\endbibitem

%b18 ###
\bibitem{xiwu11}
%
\begin{bmisc}[mr]
\bauthor{\bsnm{Xiao},~\bfnm{Han}\binits{H.}} \AND
\bauthor{\bsnm{Wu},~\bfnm{W.B.}\binits{W.B.}}
(\byear{2011}).
\bhowpublished{Simultaneous inference on sample covariances.
Available at \arxivurl{arXiv:1109.0524v1}}.
\bptok{imsref}%
\end{bmisc}
%
\endbibitem

\end{thebibliography}
\end{document}